# Stochasticity induced mixed-mode oscillations and distribution of recurrent outbreaks in an ecosystem

Susmita Sadhu[a)]
*Department of Mathematics, Georgia College & State University, Milledgeville, Georgia 31061, USA*

(Received 15 August 2016; accepted 15 February 2017; published online 7 March 2017)

The effect of stochasticity, in the form of Gaussian white noise, in a predator–prey model with two distinct time-scales is presented. A supercritical singular Hopf bifurcation yields a Type II excitability in the deterministic model. We explore the effect of stochasticity in the excitable regime, leading to dynamics that are not anticipated by its deterministic counterpart. The stochastic model admits several kinds of noise-driven mixed-mode oscillations which capture the intermediate dynamics between two cycles of population outbreaks. Depending on the strength of noise, the prey population exhibits intermediate to high-amplitude fluctuations (related to moderate or severe outbreaks respectively). We classify these fluctuations as isolated or intermittent or as clusters depending on their recurrences. We study the distribution of the random variable $N$, representing the number of small oscillations between successive spikes, as a function of the noise intensity and the distance to the Hopf bifurcation. The distribution of $N$ is "asymptotically geometric" with the corresponding parameter related to the principal eigenvalue of a substochastic Markov chain. Finally, the stochastic model is transformed into its "normal form" which is used to obtain an estimate of the probability of repeated outbreaks. *Published by AIP Publishing.*

**Large erratic changes in population densities, often referred to as outbreaks or collapses, have been documented for several species in Nature. These extreme events are typically random, with their severities varying from one event to another. This paper investigates the effect of small random perturbations in a simple bi-trophic predator–prey model, and reproduces several characteristics of the naturally occurring fluctuations in the population density of the prey, as observed in Nature. The interactions between the two species is portrayed as a system of singularly perturbed Itô stochastic differential equations. We study the effect of noise near the onset of a singular Hopf bifurcation. Remarkably, the stochastic model admits noise-driven mixed-mode oscillations that provide a realistic representation of the dynamics occurring between successive outbreaks. The random number of small oscillations between two large oscillations is assigned an integer-valued variable $N$ which can be related to the return time between two outbreaks. The distribution of $N$ is investigated to estimate the probability of successive outbreaks. This may provide useful information in preventing pest outbreaks. Since environmental fluctuations can lead to perturbations in the net growth rates of species, this study unfolds the effect of stochasticity in a stable ecosystem. Understanding population cycles of insects and small mammals with eruptive dynamics is important as recurrent pest outbreaks can lead to serious ecological and economic consequences.**

## I. INTRODUCTION

Several predator–prey interactions in which the rates of reproduction of the predators and their prey are significantly different can be modeled by systems of slow–fast differential equations.[17,26] Such systems involve a clear separation of two time-scales between the dynamics of their components. A special robust phenomenon often found in slow–fast systems is the existence of relaxation oscillators. These are limit cycles that often arise when the slow manifold is S-shaped, and the slow dynamics when restricted to the slow manifold approaches a fold point (that corresponds to a saddle-node bifurcation of the fast system), producing a jump in some other region in the phase space.[16,30] A typical relaxation oscillation cycle in a predator–prey model can be viewed as a boom and bust cycle. During such cycles, the prey population switches between two attracting branches of the slow manifold on a fast time scale, leading to its sudden outbreak or a collapse.[28] Such high amplitude, large-scale density fluctuations have been documented for small mammals as well as for many forest pest insects. Examples include voles, house mice, lemmings, gypsy moths, larch budmoths, locusts, etc.[9,12,20,25,31]

It is worth noting that the cycles of boom and bust are not necessarily periodic. In fact in some species the high amplitude fluctuations can occur multi-annually, while others show irregular outbreaks.[9,12,20] Moreover, the densities during the outbreaks or during peak phases can differ significantly from one cycle to another. It is believed that high climatic variability can limit the extent of outbreaks or disrupt the ability of predators to regulate small mammals or insects.[20] Environmental variables, such as temperature/moisture, can have differential effects on the behavior, mortality, or metabolic efficiencies of the predators and their

[a)]Electronic mail: susmita.sadhu@gcsu.edu

prey.[23] For example, low temperature decouples the ant predation of tent caterpillars in Canada, thereby increasing the probability of outbreak of tent caterpillars.[2] Unfortunately, two-dimensional deterministic predator–prey models fail to capture the intermediate dynamics between cycles of boom or bust. This paper studies such high amplitude fluctuations and the intermediate dynamics between the fluctuations by involving predator–prey interactions under varying environmental conditions. In particular, we study a stochastic version of the classical Bazykin's model[3] to reveal the role of stochasticity in generating random eruptive dynamics of the prey population.

We consider environmental variations by introducing random perturbations through Gaussian white noise in the birth rate of the prey and in the death rate of the predator. With the assumption that the per-capita growth rate of the prey is much larger than that of its predator, we portray the predator–prey interaction as a system of singularly perturbed Itô stochastic differential equations. The time diversification is measured by a small parameter $\varepsilon$. As studied in several papers,[3,22,27,28] for suitable parameter values, the deterministic model exhibits a singular Hopf bifurcation[8] at the coexistence equilibrium state. A canard explosion occurs post Hopf bifurcation, where the size of the limit cycle grows rapidly from a size of $O(\varepsilon)$ to size of $O(1)$ in an exponentially small parameter regime.[27–29]

In this paper, we are interested in the regime prior to the Hopf bifurcation, where the equilibrium state is stable, but a small external drive can "excite" the system causing an excursion away from and back to the fixed point. In fact a small perturbation in a parameter can render the equilibrium state unstable, giving birth to orbits of finite period that appear as soon as the bifurcation parameter crosses its threshold. This regime is referred to as Type II excitable regime, often encountered in neural dynamics.[19] The effect of stochasticity near the excitable regime is very interesting since it has strong effects on the dynamics which are not anticipated by its deterministic counterpart. One of the most interesting features observed in this regime are mixed-mode oscillations (MMOs). MMOs are concatenations of small and large amplitude oscillations and often occur in chemical and biological systems. In deterministic slow–fast systems, at least three variables are necessary to generate such behaviors.[7,11,15,21] However, in the presence of noise, the two-dimensional model that we consider exhibits MMOs. Noise induced MMOs have been rigorously studied for the stochastic FitzHugh–Nagumo equations.[5,24] However, such dynamics have not been yet explored in stochastic predator–prey models. The MMO cycles are ecologically significant as they qualitatively present a realistic representation of population cycles of species with eruptive dynamics. In these cycles, the two species coexist via small amplitude oscillations (SAOs) over an intermediate time scale, but the fast dynamics of the prey is randomly revealed by the burst of a large amplitude oscillation (LAO). It is to be noted that noise induced MMOs in excitable systems with one fast and two slow variables that possess solutions of canard type have been investigated in Berglund *et al.*[4]

Depending on parameters such as the growth rates, the mortality rates of the species, the noise intensity, the time scale separation, etc., we observe random fluctuations in the population densities, outbreaks of different severities, and that the frequency of recurrences of outbreaks vary. We classify them accordingly as isolated, intermittently sporadic, or repeated clusters of outbreaks in the prey population (see Figures 2 and 4). We build on the ideas by Ref. 5 to study the random variable $N$, counting the number of small oscillations between successive outbreaks, as a function of the noise intensity and the distance to the Hopf bifurcation (see Figure 13). The distribution of $N$ in our model can be associated with the relative frequency of recurrences of outbreaks. We associate the parameter corresponding to the distribution of $N$ to be related to the principal eigenvalue of a substochastic Markov chain. Finally, we transform the stochastic model into a more suitable form which we refer to as the normal form to obtain an approximate expression of the probability of repeated outbreaks.

The paper is organized as follows. In Section II, we introduce the stochastic model, rescale the parameters, and make suitable assumptions to cast the model into a framework of slow–fast system of stochastic differential equations. A preliminary review of the deterministic model is done in Section III. In Section IV, we analyze the stochastic model and discuss the distribution of the SAOs. To understand the dynamics better, the model is transformed into its normal form which is presented in Sections V and VI. We end the paper with a discussion in Section VII. All the figures in this paper have been generated through MATLAB except for Figures 1 and 8 which were done in Maple. The histograms and the time series were generated by Milstein's higher order method[18] with a step size of 0.001.

## II. THE MODEL

We consider the classical predator–prey model given by

$$\begin{cases} \dot{u} = ru\left(1 - \dfrac{u}{K}\right) - \dfrac{puv}{H+u} \\ \dot{v} = \dfrac{bpuv}{H+u} - ev - mv^2, \end{cases} \tag{1}$$

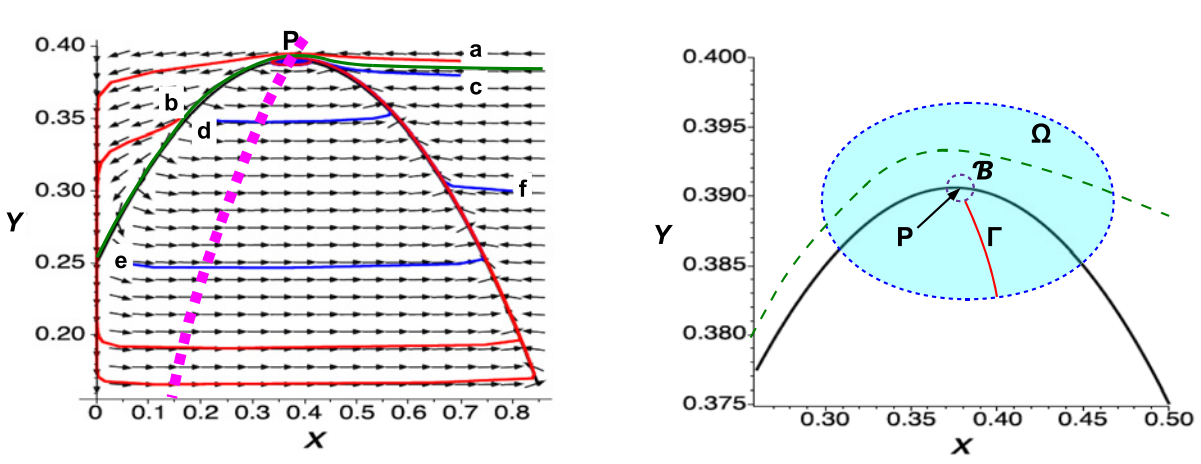

(a) (b)

FIG. 1. The parameter values used are $\beta = 0.25$, $d = 0.25$, $h = 0.91$, and $\varepsilon = 0.05$. (a) The dotted magenta curve is the $y$-nullcline and $P$ is the fixed point. Trajectories (a)–(b) and (c)–(f) start above and below the separatrix, respectively. (b) The dotted green curve represents a part of the separatrix and the black curve is the $x$-nullcline. (a) Phase portrait of system (5) exhibiting Type II excitability. (b) The bounded domain $\Omega$ containing $P$ and a part of the separatrix.

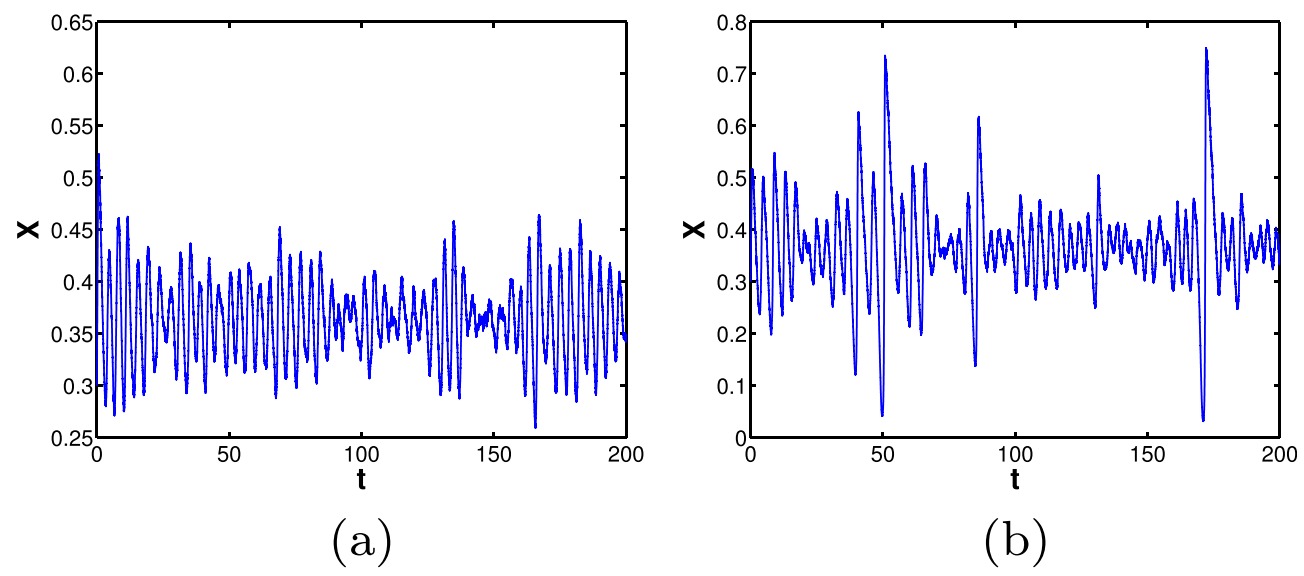

FIG. 2. Time series of the population densities of the prey corresponding to system (4) at $h = 0.875$ for weak noise intensities. The other parameter values are same as in Figure 1. (a) $\sigma_1 = \sigma_2 = 0.005$ and (b) $\sigma_1 = \sigma_2 = 0.008$.

subjected to initial conditions $u(0) > 0, v(0) > 0$. Here $u$ and $v$ represent the population densities of the prey and the predator, respectively, and $\cdot = \frac{d}{ds}$. In the absence of predators, the prey population grows logistically with an intrinsic growth rate $r$ and a carrying capacity $K$. The consumption of prey is governed by Holling type II response, with $p$ being the maximum per-capita predation rate. The parameter $H$ measures the semi-saturation constant at which the predator's predation rate is half of its maximum when $u = H$, while the parameters $b$ and $e$ represent the birth-to-consumption ratio and the per-capita mortality rate of $v$, respectively. The quantity $-mv^2$ can be interpreted as self-competition or intraspecific competition within the class of $v$. Such predator–prey interactions were first considered by Bazykin[3] and have been analyzed in detail in Refs. 22 and 28.

Incorporating environmental changes that cause random variations in the prey growth rate and in predator mortality rate, we propose the following model:

$$\begin{cases} dU = U\left[r\left(1 - \dfrac{U}{K}\right) - \dfrac{pV}{H+U} + W_1(s)\right]ds \\ dV = V\left[\dfrac{bpU}{H+U} - e - mV + W_2(s)\right]ds, \end{cases} \tag{2}$$

where $U(s)$ and $V(s)$ are two stochastic processes representing the prey and the predator population densities, respectively, and $W_1(s)$, $W_2(s)$ are two independent Gaussian white noise in the It ô sense. The auto-correlation functions of $W_1(s)$ and $W_2(s)$ are given by

$$\mathbb{E}[W_i(s)W_i(s+\tau)] = \zeta_i \delta(\tau),\ i = 1, 2,$$

where $\zeta_i$ $(i = 1, 2)$ measures the intensity of the white noise and $\delta(\tau)$ is the Dirac delta function.

System (2) can be rewritten in the form of Itô stochastic differential equation as

$$\begin{cases} dU = U\left[r\left(1 - \dfrac{U}{K}\right) - \dfrac{pV}{H+U}\right]ds + \zeta_1 U dB_1(s) \\ dV = V\left[\dfrac{bpU}{H+U} - e - mV\right]ds + \zeta_2 V dB_2(s), \end{cases} \tag{3}$$

where $B_1(s)$ and $B_2(s)$ are two independent unit Wiener processes where $dB_i = W_i ds$, $i = 1, 2$. With the following change of variables

$$t = bps, \quad X = \frac{U}{K}, \quad Y = \frac{pV}{rK}, \quad \beta = \frac{H}{K}, \quad d = \frac{e}{bp},$$
$$h = \frac{mrK}{bp^2}, \quad \varepsilon = \frac{bp}{r}, \quad \sigma_1 = \frac{\zeta_1}{\sqrt{r}}, \quad \sigma_2 = \frac{\zeta_2}{\sqrt{bp}},$$

system (3) takes the following nondimensional form:

$$\begin{cases} dX = \dfrac{X}{\varepsilon}\left[1 - X - \dfrac{Y}{\beta + X}\right]dt + \dfrac{\sigma_1}{\sqrt{\varepsilon}} X dB_1(t) \\ dY = Y\left[\dfrac{X}{\beta + X} - d - hY\right]dt + \sigma_2 Y dB_2(t). \end{cases} \tag{4}$$

Assumptions on the Parameters:

- The maximum per capita growth rate of the prey is much larger than the per-capita growth rate of its predator, i.e., $bp \ll r$, thus yielding $0 < \varepsilon \ll 1$. This is commonly observed in several predator–prey interactions in a real ecosystem.[13]
- The parameter $d$ satisfies $0 < d < 1$ which is a default assumption; otherwise, the predator would die out faster than they could reproduce even at their maximum reproduction rate.
- The parameter $\beta$ represents the dimensionless semi-saturation constant measured against the prey's carrying capacity. We will assume that the predator is efficient, and hence will reach half of its maximum predation rate before the prey population reaches its carrying capacity which leads to $0 < \beta < 1$.

**Remark II.1.** *(a) Under the first assumption, we note that system (4) represents a stochastic slow–fast system with $X$ being the fast variable and $Y$ being the slow one. The scaling of the noise intensities takes into account the difference between the time-scales. We assume that the diffusive nature of the Brownian motion causes paths to spread like $\sigma_i\sqrt{t}$, $i = 1, 2$.*

*(b) The quantities $\sigma_1^2$ and $\sigma_2^2$ measure the ratios between the rate of diffusion squared and the speed of drift for the fast and slow variables, respectively. More specifically in this context, $\sigma_1^2$ $(\sigma_2^2)$ measures the ratio between the square of the noise intensity associated with the random variation in the growth rate of the prey (death rate of the predator) and the intrinsic growth rate of the prey (predator) in the absence of environmental fluctuation.*

*(c) The parameter $h$ will be interpreted as the strength of the intraspecific competition within $Y$. This will be treated as our bifurcating parameter in the analysis.*

*(d) Similar scaling variables were first considered in a deterministic food chain model by Deng[10] and have been also incorporated in* Refs. 27–29.

## III. A REVIEW OF THE DETERMINISTIC MODEL

In this section, we consider the deterministic counterpart of system (4), namely

$$\begin{cases} \varepsilon\dot{x} = x\left(1 - x - \dfrac{y}{\beta + x}\right) \\ \dot{y} = y\left(\dfrac{x}{\beta + x} - d - hy\right), \end{cases} \tag{5}$$

where $x$ and $y$ represent the prey and the predator population densities, respectively. The slow–fast system (5) admits a slow manifold given by $\{x=0\}\cup\{y=(1-x)(\beta+x)\}$. Depending on the parameters, the nullclines $\{y=(1-x)(\beta+x)\}$ and $\{hy=\frac{x}{\beta+x}-d\}$ can have one, two, or three intersections. In this analysis, we will keep $\beta$ and $d$ fixed and allow $h$ to be the varying parameter. We consider the situation when the nullclines intersect at a unique stationary point $P=\big(\alpha,(1-\alpha)(\beta+\alpha)\big)$, where $\alpha$ satisfies the relation

$$\frac{\alpha}{\beta+\alpha}-d=h(1-\alpha)(\beta+\alpha). \tag{6}$$

The Jacobian matrix of the vector field at $P$ is

$$J=\begin{pmatrix}\frac{\alpha}{\varepsilon}\left(\frac{1-\beta-2\alpha}{\beta+\alpha}\right) & -\frac{\alpha}{\varepsilon(\beta+\alpha)}\\ \frac{\beta(1-\alpha)}{\beta+\alpha} & -h(1-\alpha)(\beta+\alpha)\end{pmatrix}. \tag{7}$$

It has trace

$$Tr(J)=\frac{\alpha^\star-\alpha}{\varepsilon(\beta+\alpha)(\beta+\alpha^\star)}\left[2\alpha\alpha^\star+2(\alpha+\alpha^\star)\beta-\beta+\varepsilon\beta+\beta^2\right] \tag{8}$$

and determinant

$$\det(J)=\frac{\alpha(1-\alpha)}{\varepsilon}\left[\frac{\beta}{(\beta+\alpha)^2}-h(1-\beta-2\alpha)\right], \tag{9}$$

where

$$\alpha^\star=\frac{(1-\beta)-\varepsilon(1-d)+\sqrt{\big((1-\beta)-\varepsilon(1-d)\big)^2+8d\beta\varepsilon}}{4}.$$

We note that $P=\big(\alpha^\star,(1-\alpha^\star)(\beta+\alpha^\star)\big)$ at $h=h^\star(\varepsilon)$, and as $\varepsilon\to0$, $P\to\Big((1-\beta)/2,(1+\beta)^2/4\Big)$ and $h^\star(\varepsilon)\to\tilde{h}$, where $h^\star$ and $\tilde{h}$ are defined by

$$h^\star=\frac{\alpha^\star(1-\beta-2\alpha^\star)}{\varepsilon(1-\alpha^\star)(\beta+\alpha^\star)^2}, \tag{10}$$

$$\tilde{h}=\frac{4(1-\beta-d(1+\beta))}{(1+\beta)^3}. \tag{11}$$

At $\alpha=\alpha^\star$, $Tr(J)|_{\alpha^\star}=0$ and $\det(J)|_{\alpha^\star}>0$ (see Appendix A), provided that

$$\frac{(1-d)^2}{\beta}<\frac{1}{\varepsilon},\quad d\beta<(1-d)(1-\beta) \tag{12}$$

hold. Hence if (12) holds and $h=h^\star$, $J$ admits a pair of conjugate imaginary eigenvalues when $\alpha=\alpha^\star$. The real parts of the eigenvalues of $J$ have the order $(\alpha-\alpha^\star)/\varepsilon$ near $\alpha^\star$. The system undergoes a singular Hopf bifurcation at $\alpha=\alpha^\star$.

We are interested in the regime when $\alpha-\alpha^\star$ is small and positive. For $\alpha>\alpha^\star$, we can show that $Tr(J)<0$ and $\det(J)>0$ provided that $\varepsilon>0$ is sufficiently small and

$$d<\frac{1-\beta}{1+\beta} \tag{13}$$

holds (see Appendix A). Hence under these conditions, $P$ is a stable stationary point for all $\alpha>\alpha^\star$, corresponding to the coexistence state. In fact the fixed point $P$ corresponds to a focus when $0<\alpha^\star-\alpha=O(\sqrt{\varepsilon})$.

The regime where $\alpha-\alpha^\star$ is small and positive is said to be excitable because a small change in a parameter can make $P$ unstable, while creating a stable periodic orbit. The Hopf bifurcation involved here yields Type II excitability,[6] where orbits of finite but rapidly growing period appear (also known as a canard explosion) as soon as the bifurcation parameter $h$ crosses its threshold value $h^\star$. However, more importantly, even before the bifurcation, trajectories that closely miss the stable equilibrium $P$ make a large excursion before reaching $P$. As seen in Figure 1(a), a small perturbation of the initial condition can cause the system to make a large excursion corresponding to a spike, leading to huge changes in the population densities of the species, before the system reaches $P$. In fact, there exists a separatrix (represented by the green curve) that separates the long excursions and the short amplitude oscillations. The separatrix can be thought of as a trajectory closely tracking the unstable branch of the slow manifold $\{y=(1-x)(\beta+x)\}$.

In the presence of noise, the system generates sample paths making excursions away from $P$. A sample path can reach and cross the separatrix, causing the system to spike. The shape and duration of spikes will be close to their deterministic value, but the time between the spikes will be random.

## IV. ANALYSIS OF THE STOCHASTIC MODEL

We will work in the excitable regime, namely, all those $\alpha>\alpha^\star$ that are $O(\sqrt{\varepsilon})$ distance away from $\alpha^\star$. For $\alpha>(1-\beta)/2$, the slow manifold of the deterministic system is uniformly asymptotically stable. Hence Fenichel's theorem[14] ensures the existence of an invariant manifold $\mathcal{M}_\varepsilon$. Sample paths stay in an appropriate constructed neighborhood of $\mathcal{M}_\varepsilon$ for exponential long time spans until the system is driven towards a bifurcation point, where the slow manifold ceases to be attracting. Detailed analysis of the dynamics with precise estimates near the slow manifold has been done in Ref. 6. In this paper, we consider the regime near the singular Hopf bifurcation. A similar analysis has been previously done for stochastic FitzHugh–Nagumo equations in Ref. 5. For chosen parameter values, we observe LAOs separated by random number of SAOs. These oscillations are referred to as MMOs.[7,11,15,21] The LAOs in our context signify cycles of boom and bust, while the SAOs represent small fluctuations in the population densities. The MMOs provide us with a more realistic situation, taking into account the intermediate dynamics between sudden outbreaks or crashes which may occur due to environmental variations.

The excitable regime is ecologically significant. In this regime in the absence of noise, the deterministic model exhibits a stable stationary point, corresponding to the coexistence state of the two species. However, when subjected to small environmental variations, the stationary point loses its stability and the species exhibit small fluctuations in their densities, interspersed with large fluctuations at times as shown in Figure 2.

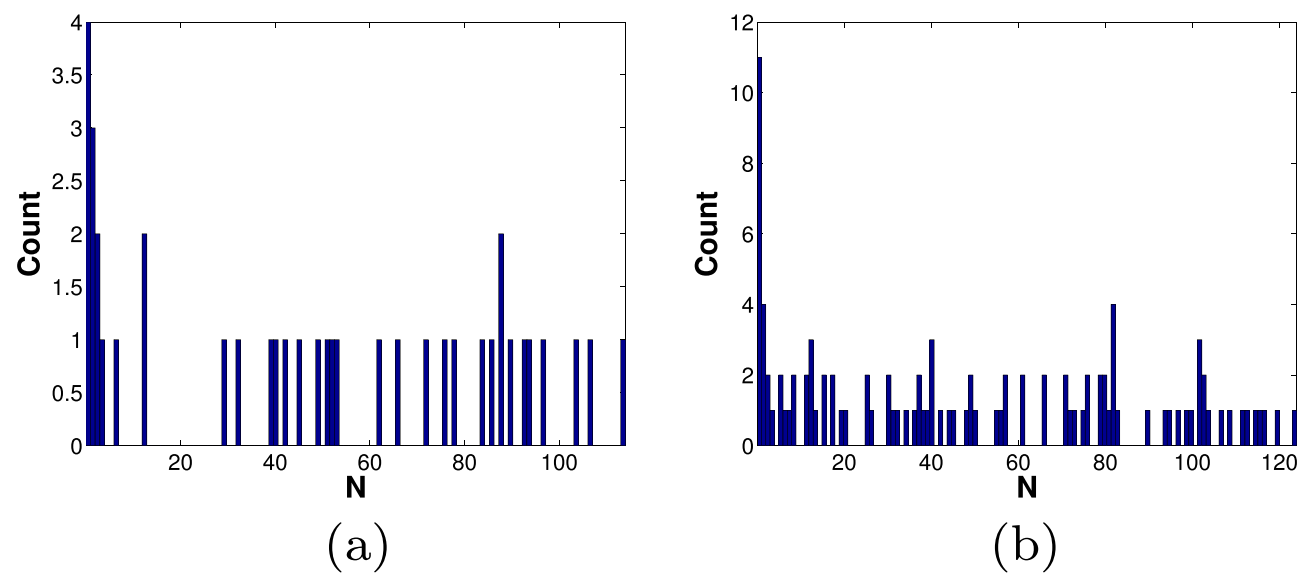

FIG. 3. Histograms of the distribution of the number of SAOs (represented by $N$) between two spikes of system (4) performed over 200 sample paths over the interval [0, 500]. The other parameter values are same as in Figure 2. Here we chose $m_t = 0.68$, $a_t = 0.06$, and $a_t' = 0.3$ (see Remark IV.1). (a) $(\sigma_1 = \sigma_2) = (0.007, 0.005)$, $\mathbb{E}(N) = 47.97$, and std = 38.15. (b) $(\sigma_1 = \sigma_2) = (0.006, 0.006)$, $\mathbb{E}(N) = 48.74$, and std = 38.84.

For weak noise intensities, the prey population fluctuates around a stationary state performing a large number of SAOs before making an LAO. The distribution of the number of SAOs between two spikes for such noise intensities is shown in Figure 3 (STD stands for standard deviation). The probability of frequent outbreaks is very low in this situation.

**Remark IV.1.** *For numerical purposes, we separate SAOs from LAOs by detecting the amplitude of the oscillations. Some of the SAOs are indistinguishable from the random fluctuations, and so we set a threshold on the amplitudes of oscillations. All such oscillations whose amplitudes exceed a threshold value* $a_t$ *with their maxima below a threshold value* $m_t$ *are regarded as SAOs. Oscillations with amplitudes greater than* $a_t'$ *for some* $a_t' > a_t$ *and maxima greater than* $m_t$ *are designated as LAOs.*

Figure 4 shows the effect of intermediate and/or strong noise that interplays intricately with the deterministic dynamics leading to qualitatively distinct behaviors. The transition from excitability to LAOs (relaxation oscillations) is highly nontrivial, capturing some of the dynamics as observed in Nature. For example, a small perturbation in the birth rate of the prey may lead to a sudden outbreak in its population, which then gradually collapses due to overpopulation or starvation. The population may then stay at a dormant mode until a favorable environmental condition triggers another outbreak (see (a) and (b) in Figure 4). Here outbreaks and collapses correspond to the highest and the lowest values of the spikes, respectively. At an intermediate noise intensity (see (c) in Figure 4), the frequencies of the outbreaks increase, thus relating to outbreaks of insects, such as gypsy moths, spruce budworms, larch budmoths, etc.[1,12] Higher noise intensities (see (d) in Figure 4), on the other hand, can lead to more frequent and sometimes consecutive outbreaks as recorded for Australian plague locust.[9,31] The population densities during outbreaks, which signify the severity of an outbreak, vary from cycle to cycle. Some outbreaks can be regarded as moderate, while others as severe, depending on the environmental perturbations (see (b) and (c) in Figure 4).

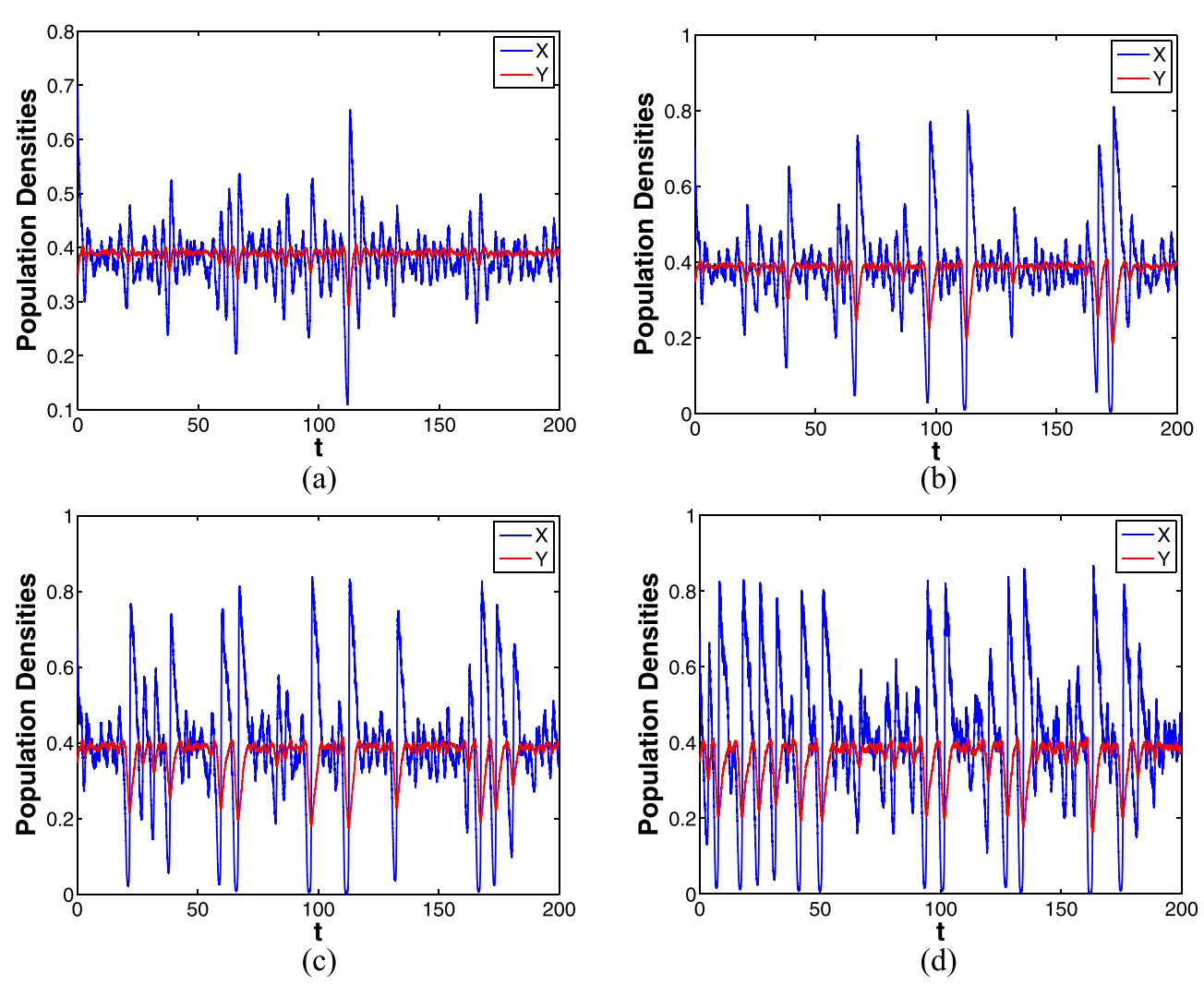

FIG. 4. Time series of sample paths of the population densities corresponding to system (4) at $h = 0.9$. The other parameter values are same as in Figure 1. (a) An isolated outbreak in the prey population when $\sigma_1 = \sigma_2 = 0.01$. (b) Intermittent sporadic outbreaks when $\sigma_1 = \sigma_2 = 0.014$. (c) An isolated cluster of outbreaks when $\sigma_1 = \sigma_2 = 0.02$. (d) Repeated clusters of outbreaks when $\sigma_1 = \sigma_2 = 0.03$.

### A. The distribution of the SAOs

Similar to the work in Ref. 5, we can analyze the general properties of the integer-valued random number $N$, counting the number of SAOs between two consecutive LAOs that system (4) exhibits. To this end, we fix a bounded set $\Omega \subset \mathbb{R}^2$, with smooth boundary $\partial\Omega$, containing the stationary point $P$ and a piece of the separatrix. Let $\mathcal{B}$ be a ball of radius $r > 0$, centered at $P$, with $r$ being small. Let $\Gamma$ be a smooth curve drawn from a point in $\mathcal{B}$ to a point in $\partial\Omega$ (see Figure 1(b)). The curve $\Gamma$ is parametrized by a variable $r$ proportional to arc length. The parametrization is then extended to all $\Omega\backslash\mathcal{B}$ by choosing a diffeomorphism $T : [0, 1/2] \times \mathbb{S}^1 \to \Omega$ by $(r, \phi) \mapsto (x, y)$, where $T^{-1}(\Gamma) = \{\phi = 0\}$, $T^{-1}(\partial\Omega) = \{r = 0\}$ and $T^{-1}(\partial\mathcal{B}) = \{r = 1/2\}$. We also arrange that $\dot{\phi} > 0$ near $P$ for the deterministic flow.

We consider the process $(r_t, \phi_t)_t$. Given an initial condition $(r_0, 0) \in T^{-1}(\Gamma)$ and an integer $M \geq 1$, we define the stopping time

$$\tau = \inf\{t > 0 : \phi_t \in \{2\pi, -2M\pi\} \text{ or } r_t \in \{0, 1/2\}\}.$$

As interpreted in Ref. 5, the case $r_\tau = 0$ corresponds to the sample path $(X_t, Y_t)$ leaving $\Omega$, thus giving rise to an LAO. We set $N = 0$ in this situation. The case $r_\tau = 1/2$ corresponds to the situation when the sample path enters $\mathcal{B}$ and reaches a dormant state and stays close to $P$. In this case, we wait until the state leaves $\mathcal{B}$ and either hits $\Gamma$ or leaves $\Omega$. When $\phi_\tau = 2\pi$ and $r_\tau \in (0, 1/2)$, the sample path has returned to $\Gamma$ after performing a complete revolution around $P$. In this situation, $N \geq 1$. Finally, the case $\phi_\tau = -2M\pi$ and $r_t \in (0, 1/2)$ represents the (unlikely) event that the sample path winds $M$ times around $P$ in the wrong direction. For simplicity, in this situation we consider $N = 1$.

As long as $r_\tau \in (0, 1/2)$, the procedure is repeated, incrementing $N$ at each iteration, thus yielding a sequence $(R_0, R_1, \ldots R_{N-1})$ of random variables. The sequence corresponds to the successive intersections of the path with $\Gamma$, separated by rotations around $P$, until it first exits from $\Omega$ forming a substochastic Markov chain on $E = (0, 1/2)$. Suppose that $R \in E$ and $A \subseteq E$ is a Borel set. Then the kernel $K(R, A)$ is defined as

$$K(R, A) = \mathbb{P}\{\phi_\tau \in \{2\pi, -2M\pi\}, r_\tau \in A \,|\, \phi_0 = 0, r_0 = R\}.$$

The Markov chain is substochastic because $K(R, E) < 1$, due to the positive probability of sample paths leaving $\Omega$. By adding a cemetery state $\Delta$ to $E$ and setting $K(R, \Delta) = 1 - K(R, E)$, $K(\Delta, \Delta) = 1$, we can make the Markov chain stochastic. The number of SAOs is given by

$$N = \inf\{n \geq 0 : R_n = \Delta\} \in \mathbb{N} \cup \{0\} \cup \{\infty\},$$

where we set $\inf \emptyset = \infty$. A suitable extension of the Perron–Frobenius theorem implies the existence of the principal eigenvalue $\lambda_0$ of $K$. Theorem 3.2 in Ref. 5 proves that for $\sigma_1$, $\sigma_2 > 0$ in (4) and for any initial distribution $\mu_0$ of $R_0$ on the curve $\Gamma$, the kernel $K$ admits a quasi-stationary distribution $\pi_0$ satisfying $\pi_0 K = \lambda_0 \pi_0$, where $\lambda_0 < 1$. The random variable $N$ is almost surely finite, i.e., $\lim_{n\to\infty} \mathbb{P}^{\mu_0}\{N > n\} = 0$; the distribution of $N$ is "asymptotically geometric" (see Figure 5), i.e.,

$$\lim_{n\to\infty} \mathbb{P}^{\mu_0}\{N = n + 1 \,|\, N > n\} = 1 - \lambda_0;$$

and the probability-generating function $\mathbb{E}^{\mu_0}\{\theta^N\} < \infty$ for $\theta < 1/\lambda_0$ which implies that $N$ has finite moments (see Figure 6). The probability generating function $\mathbb{E}^{\mu_0}\{\theta^N\}$ has a simple pole at $\theta = 1/\lambda_0$. For practical purposes, we can estimate the value of $\lambda_0$ by detecting when the derivative of the probability generating function exceeds a given threshold (see Figures 6 and 7). In Figure 7, we plot the numerically simulated probability mass function of $N$ of system (4) (blue curve) and the probability mass function of a geometric distribution with parameter $1 - \lambda_0$ (red curve).

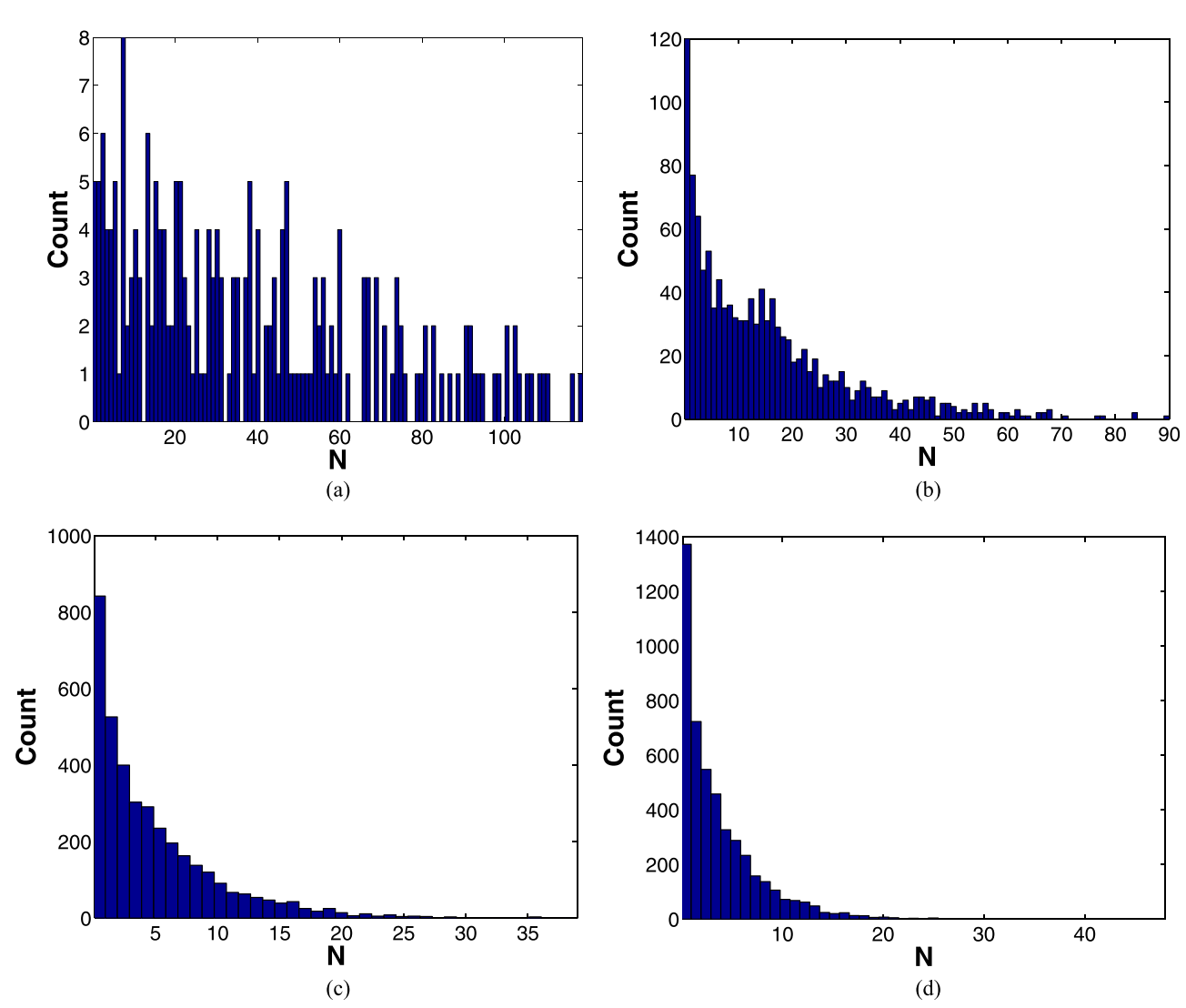

FIG. 5. Histograms of the distribution of the $N$ corresponding to system (4) generated over 200 sample paths with $\sigma_1 = \sigma_2 = 0.015$. The other parameter values are as in Figure 3. (a) $h = 0.92$, $\mathbb{E}(N) = 38.71$, and std $= 30.81$. (b) $h = 0.9$, $\mathbb{E}(N) = 15.1$, and std $= 15.14$. (c) $h = 0.88$, $\mathbb{E}(N) = 4.5$, and std $= 5.18$. (d) $h = 0.875$, $\mathbb{E}(N) = 3.3$, and std $= 4.06$.

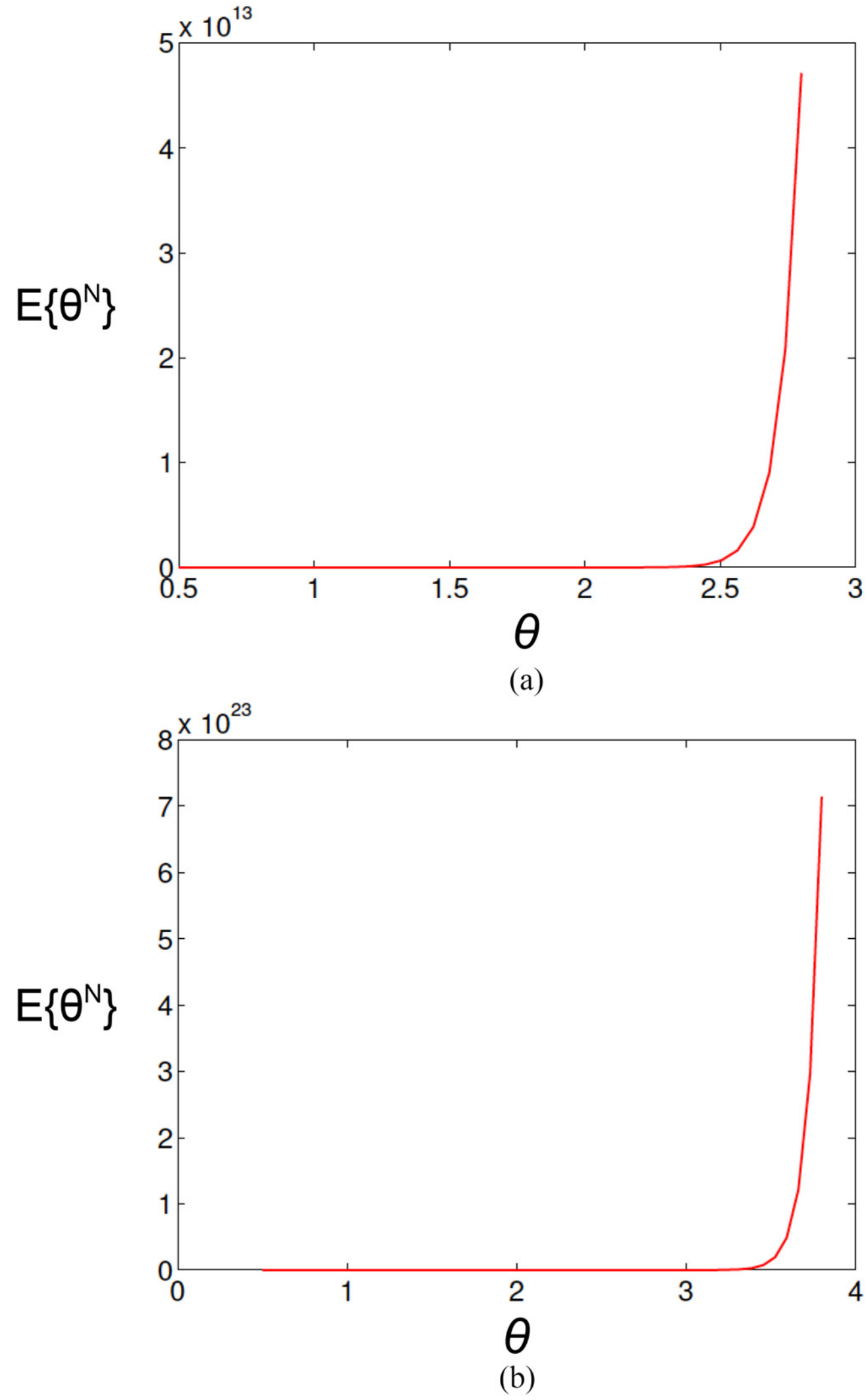

FIG. 6. Probability-generating functions of $N$ of system (4) with other parameter values as in Figure 5. (a) $h = 0.88$ and (b) $h = 0.875$.

## V. NORMAL FORM OF THE DETERMINISTIC MODEL

Following the ideas of Berglund and Landon,[5] we will transform system (5) into a suitable form which we refer to as the "normal form." The normal form allows us to obtain more quantitative results. With this transformation, we introduce a new coordinate system $(l, z)$ such that the corresponding dynamical system admits a separatrix in the lower half plane close to the horizontal axis $z = 0$. The separatrix delimits the quiescent and the spiking regimes. The new system makes it easier to separately analyze the dynamics near the separatrix. We are interested in the parameter regime near the singular Hopf bifurcation point. We assign the distance to the Hopf point, denoted by $\mu$, as the primary bifurcation parameter. At the threshold value of $\mu$, i.e., when $\mu = 0$, the normal form admits a first integral when $\varepsilon = 0$. For $\mu \neq 0$, in the vicinity of the fixed point, the dynamics stay close to the level curves of the first integral.

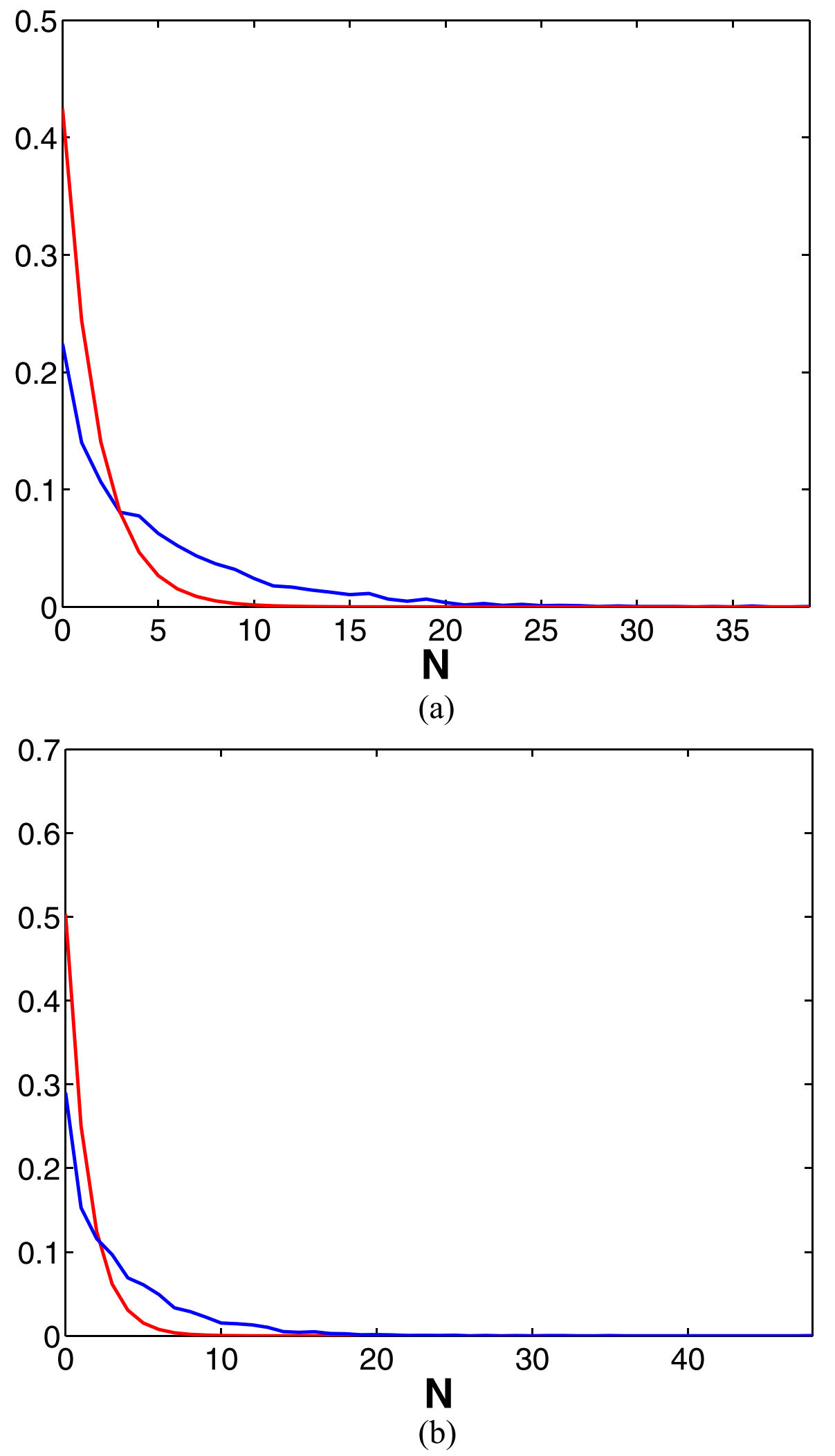

FIG. 7. The probability mass functions of $N$. (a) $h = 0.88$, $\lambda_0 = 0.5752$. (b) $h = 0.875$, $\lambda_0 = 0.4969$.

**Proposition V.1.** *Under suitable transformations as provided in Appendix B and for $\varepsilon > 0$, system (5) is orbitally equivalent to the following system of equations:*

$$\begin{cases} \dot{l} = \dfrac{1}{2} - z + A_0(l,z)\sqrt{\varepsilon} \\ \dot{z} = \mu + 2\alpha^\star \gamma l z + A_1(l,z)\sqrt{\varepsilon} + A_2(l,z)\varepsilon \\ \qquad + A_3(l,z)\varepsilon^{3/2} + A_4(l,z)\varepsilon^2, \end{cases} \tag{14}$$

*where*

$$\mu = \frac{1-\beta-3\alpha^\star}{(1-d-h(1-\alpha^\star)(\beta+\alpha^\star))}\frac{\delta}{\sqrt{\varepsilon}}, \tag{15}$$

$$\delta = \alpha^\star - d(\beta+\alpha^\star) - h(1-\alpha^\star)(\beta+\alpha^\star)^2, \tag{16}$$

$$\gamma = (1-\alpha^\star)(\beta+\alpha^\star)\big(1-d-h(1-\alpha^\star)(\beta+\alpha^\star)\big), \tag{17}$$

*and $A_i$, $i = 0\ldots 4$ are bi-variate polynomials in $l$ and $z$ defined by (A4)–(A8) in Appendix B.*

**Remark V.1.** *The parameter $\mu$ as defined by (15) is associated with the distance to the Hopf bifurcation point (see Equation (A3)). From the definition of $\alpha^\star$ and for $h^\star < h < \tilde{h}$, it follows that $\mu > 0$ (follows from equation (16)), provided that (12) and (13) hold.*

We will first consider some special cases below.

Case (I): Ignoring the terms in $\varepsilon$, system (14) reduces to

$$\begin{cases} \dot{l} = \dfrac{1}{2} - z \\ \dot{z} = \mu + 2\alpha^\star \gamma l z. \end{cases} \tag{18}$$

(a) $\mu = 0$: In this situation, system (18) admits a first integral given by

$$Q = z e^{-2\alpha^\star \gamma l^2 - 2z + 1}. \tag{19}$$

Note that $Q = 1/2$ corresponds to the fixed point $(0, 1/2)$ while $0 < Q < 1/2$ to periodic orbits. The level curve $Q = 0$ represents the separatrix $z = 0$. Trajectories below $z = 0$ diverge to $-\infty$ in time (see Figure 8(a)).

(b) $\mu \neq 0$: In this situation, the stationary point moves to $P = \left(-\frac{\mu}{\alpha^\star\gamma}, \frac{1}{2}\right)$. When $\mu > 0$, $P$ is a stable sink, while if $\mu < 0$, $P$ corresponds to an unstable focus. When $\mu > 0$, the separatrix is deformed and lies in the lower half plane (see Figure 8(b)). Orbits below the separatrix diverge to $-\infty$.

Case (II): $0 < \varepsilon \ll 1$, $\mu > 0$: The local dynamics stays same as in Case (I)(b) above. The stationary point $P$ is a stable sink. The global dynamics is topologically equivalent to system (5). Orbits below the separatrix are no longer unbounded, but make a large excursion in the plane before converging to $P$. This large excursion corresponds to a large amplitude oscillation. Orbits above the separatrix converge to $P$ in an oscillatory manner (see Figure 9).

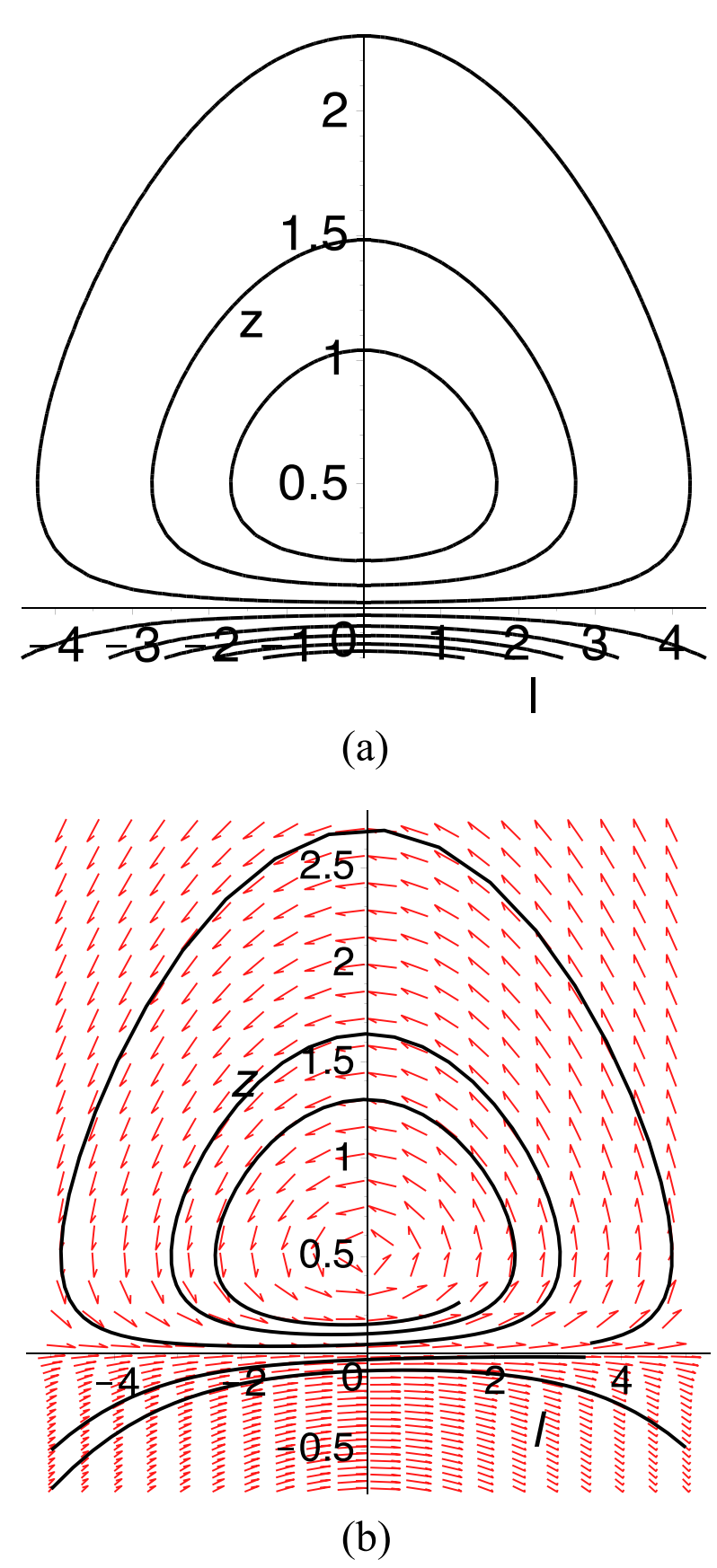

FIG. 8. Integral curves of system (18). (a) $\mu = 0$ and (b) $\mu = 0.0062$.

## VI. NORMAL FORM OF THE STOCHASTIC MODEL

In this section, we apply the same series of transformation as done in the deterministic model to the stochastic model (4). We apply Itô's formula in the derivation. The new coordinate system with variables ($L$, $Z$) will help us to describe the dynamics in a neighborhood of the separatrix. The variables $L$ and $Z$ are the stochastic counterparts of $l$ and $z$, respectively. Whether the system performs a large amplitude oscillation depends on the dynamics near the separatix which is close to the horizontal line $Z=0$. The new system is given by Proposition IX.1 in Appendix C.

The Itô-to-Stratonovich correction terms (that appear in Proposition IX.1) are multiples of $-\hat{\sigma_1}^2$, defined by (A10) and occur in the expressions of $\hat{\mu}$, $B_1$, $B_2$, and $B_3$ defined by (A11), (A14), (A15), and (A16), respectively. The parameter $\hat{\mu}$ can be positive or negative depending on the value of $\delta$ and the noise intensity. For small values of $\varepsilon$ and $\hat{\mu}$ and for low noise intensities, the sample paths corresponding to system (A9) stay close to the level curves of the first integral $Q$ and thus oscillate near the stationary point. At intermediate or higher noise intensities, the system exhibits MMOs, namely, short amplitude oscillations interspersed with long excursions as shown in Figure 10. The trajectory spends most of its time performing SAOs by oscillating around the deterministic fixed point until it is pushed off to the other side of the separatrix by the noise.

To this end, we consider the dynamics near the separatrix. Consider a sample path of (A9) starting with an initial condition ($L_0$, $Z_0$) where $L_0=-P$ for some $P>0$ (which will be determined later) and $Z_0$ small. Ignoring the terms in $\varepsilon$ in (A9), we consider the following system:

$$\begin{cases} dL = \frac{1}{2}dt' \\ dZ = (\hat{\mu} + 2\alpha^\star\gamma LZ)dt' + \sqrt{\beta+\alpha^\star}(-\gamma\alpha^\star\hat{\sigma}_1 dB_1(t') \\ \qquad +(1-\alpha^\star)(\beta+\alpha^\star)\hat{\sigma}_2 dB_2(t')). \end{cases}$$

Note that as long as $Z_{t'}$ remains small, we can approximate $L_{t'}$ in the mean by $L_0+t'/2$, and so $Z_{t'}$ will be close to the solution of

$$dZ_{t'} = (\hat{\mu} + \gamma\alpha^\star t Z_{t'})dt' + \sqrt{\beta+\alpha^\star}\big(-\gamma\alpha^\star\hat{\sigma}_1 dB_1(t') + (1-\alpha^\star)(\beta+\alpha^\star)\hat{\sigma}_2 dB_2(t')\big). \quad (20)$$

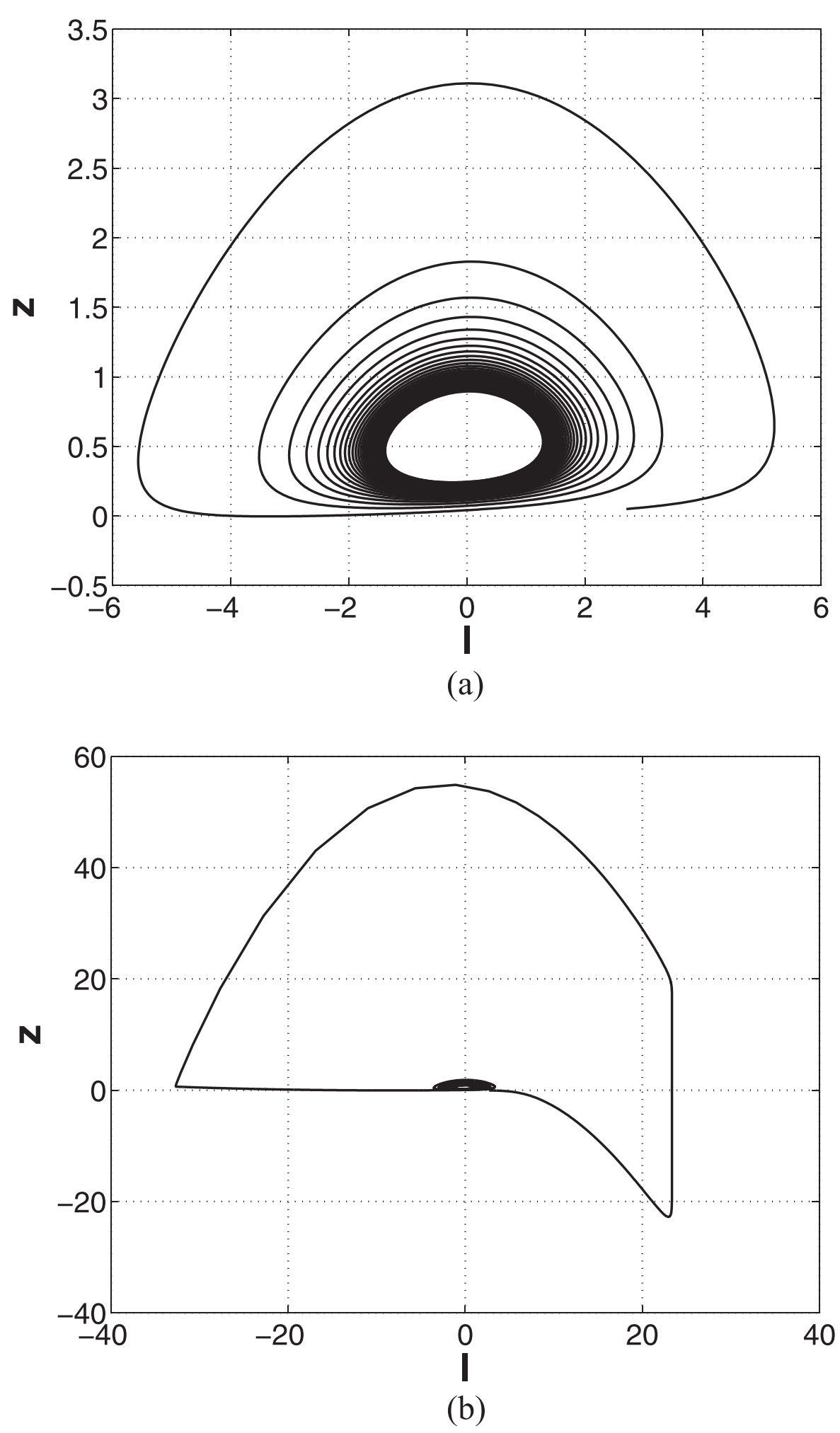

FIG. 9. Dynamics of system (14) for $\mu=0.0017$, $\varepsilon=0.01$, $\beta=0.25$, and $d=0.25$. IC: initial conditions. (a) IC: (2.7, $-0.05$) and (b) IC: (2.7, $-0.05$).

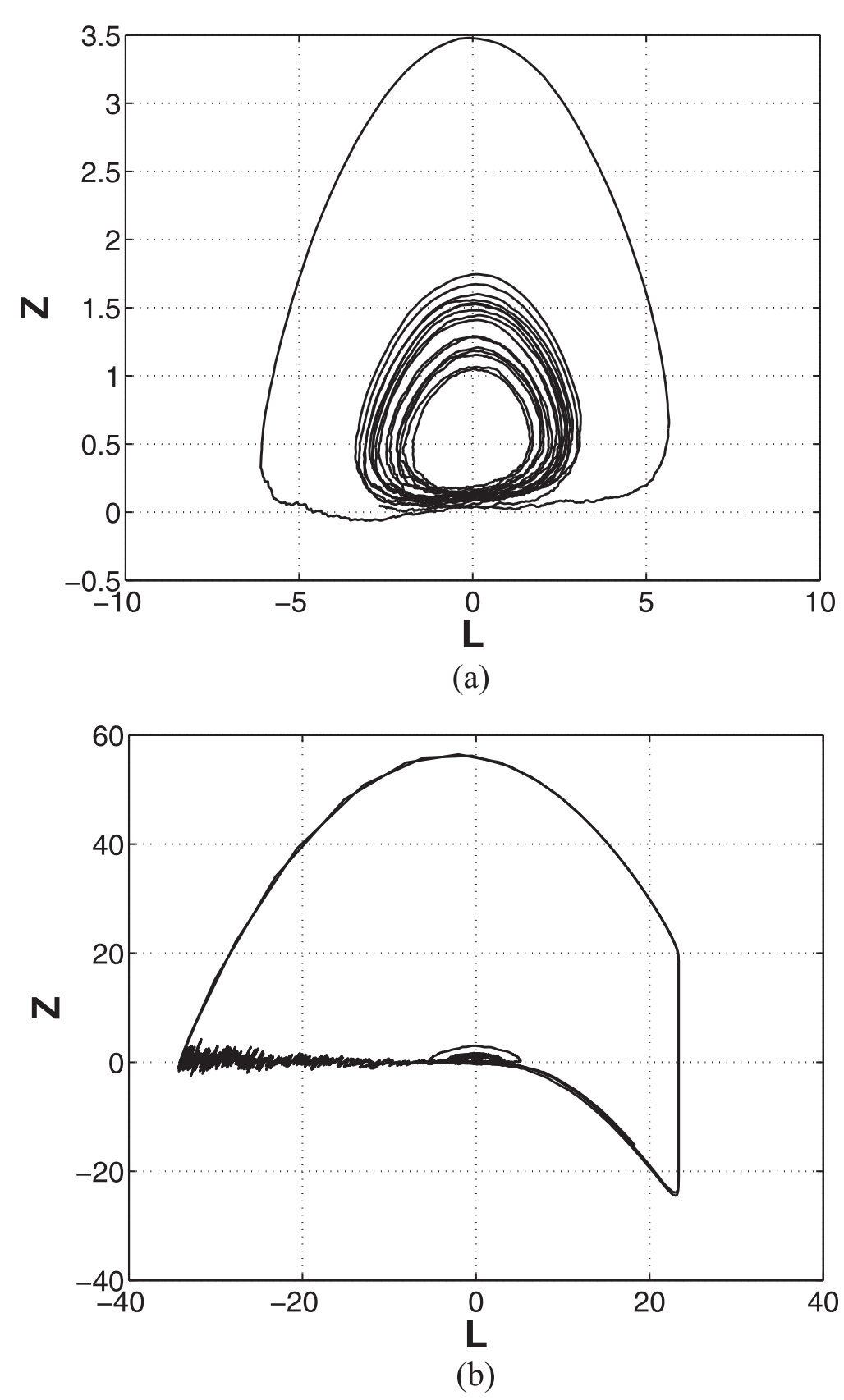

FIG. 10. Sample paths of system (A9) for parameter values as in Figure 9. (a) $\sigma_1=\sigma_2=0.0005$. (b) $\sigma_1=\sigma_2=0.005$.

We can find an estimate of the probability of repeated spiking using the proposition below.

**Proposition VI.1.** *Let* $2\gamma\alpha^\star P^2 = \log|(c_0\hat{\mu})^{-a}|$ *for some* $a \in (0,1)$ *and* $c_0 > 0$. *Then for any* $H > 0$

$$\mathbb{P}\{Z_T \le -H\} \le \Phi\big(-\kappa\big(1 + O(H + Z_0)\hat{\mu}^{a-1}\big)\big),$$

*where*

$$\kappa = \frac{\hat{\mu}\pi^{1/4}}{(\gamma\alpha^\star)^{1/4}\sqrt{\left(\frac{\gamma\alpha^\star\hat{\sigma}_1{}^2}{2} + (1-\alpha^\star)^2(\beta+\alpha^\star)^2\hat{\sigma}_2{}^2\right)(\beta+\alpha^\star)}}$$

*and*

$$\Phi(x) = \frac{1}{\sqrt{2\pi}}\int_{-\infty}^{x} e^{-\frac{u^2}{2}}du,$$

*is the cumulative probability distribution function of the standard normal law.*

The proof is provided in Appendix D.

For large values of $H$, the system performs no complete SAOs before spiking. Hence by choosing $a$ large enough, and by considering large values of $H$, we note that $\mathbb{P}\{Z_T \le -H\}$ can be estimated by $\Phi(-\kappa)$. Therefore, the probability of repeated spiking $\mathbb{P}^{\mu_0}\{N = 0\} \approx \Phi(-\kappa)$. If the noise intensity and/or distance to the Hopf bifurcation be such that $\hat{\mu} < 0$, then $\kappa < 0$ leads to a higher probability of repeated spiking (see Figure 11). This is expected, since the closer we are to the Hopf bifurcation, the more excitable the system is. With an increasing noise intensity, the recurrences of spikes increase. The contour plot of repeated outbreaks as a function of $\sigma = \sqrt{\sigma_1^2 + \sigma_2^2}$ and $\mu$ is presented in Figure 13(b).

**Remark VI.1.** If $\hat{\mu} \gg \hat{\sigma}$, *where* $\hat{\sigma} = \sqrt{\hat{\sigma}_1{}^2 + \hat{\sigma}_1{}^2}$, *which in original variables translates to* $\sigma_1^2 + \sigma_2^2 \ll \sqrt{\varepsilon}(1-\alpha^\star)^2(\beta+\alpha^\star)^2\delta^2$, *then* $\kappa \gg 0$, *and hence we can expect that* $\mathbb{P}^{\mu_0}\{N = 0\} \approx 0$. *On the other hand, if* $O(\hat{\mu}) = O(\hat{\sigma})$, *then* $\mathbb{P}^{\mu_0}\{N = 0\} \approx 1/2$, *while if* $\hat{\sigma}$ *be such that* $\hat{\mu} < 0$, *then* $\mathbb{P}^{\mu_0}\{N = 0\} > 1/2$ *(see Figure 13(b))*.

## VII. DISCUSSION AND SUMMARY

A randomized Bazykin's model is used to analyze the interactions between the predators and their prey in an ecosystem. We model the random variations in the birth rate of the prey and the death rate of the predator through Gaussian white noise. Considering predator–prey interactions in a bi-trophic ecosystem, where the prey exhibits a faster dynamics than that of its predator, we cast our model as a system of singularly perturbed Itô SDEs. We focus our analysis in the excitable regime, where the coexistence equilibrium state is stable in the absence of noise, but the system is sensitive to small perturbation.

The effect of stochasticity in the excitable regime is mathematically as well as ecologically interesting. In the excitable regime, depending on the strength of noise, the system exhibits MMOs. To the best of our knowledge, noise induced MMOs have not yet been explored in two-dimensional stochastic predator–prey models. The time series of the prey population for the parameter values in this regime give a good qualitative representation of population dynamics of species exhibiting high amplitude fluctuations in their densities. We study the random number of SAOs, denoted by $N$, in between two spikes. Here we relate the SAOs to small fluctuations in the population densities of the prey that occur in between two outbreaks. We relate the random variable $N$ to a substochastic continuous-space Markov chain. The distribution of $N$ is asymptotically geometric with parameter $1 - \lambda_0$, where $\lambda_0$ is related to the principal eigenvalue of the Markov chain. For large values of $N$, the distribution closely follows the exponential law with parameter $1 - \lambda_0$ as shown in Figure 7. The distribution of $N$ can be used to obtain the distribution of return time of the outbreaks.

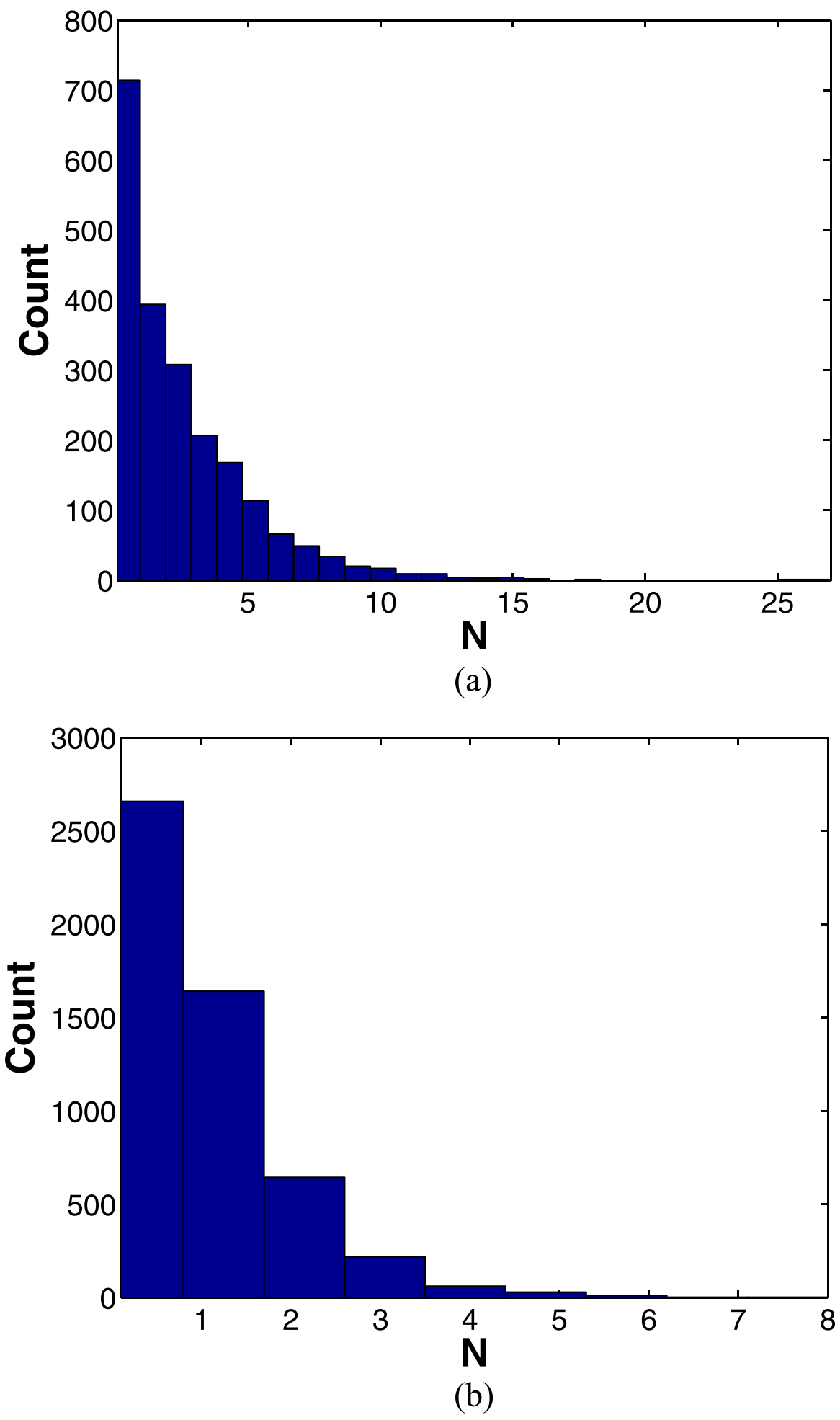

FIG. 11. Histograms of distribution of $N$ for system (A9) computed for 200 sample paths over [0, 2000] with other parameter values as in Figure 5. Here $m_t = 6$, $a_t = 0.1$, and $a'_t = 3$. (a) $\hat{\sigma}_1 = \hat{\sigma}_2 = -0.2966$, $\hat{\mu} = 0.0036$, and $\Phi(-\kappa) \approx 0.4600$. (b) $\hat{\sigma}_1 = \hat{\sigma}_2 = -0.9802$, $\hat{\mu} = -0.035$, and $\Phi(-\kappa) \approx 0.6146$.

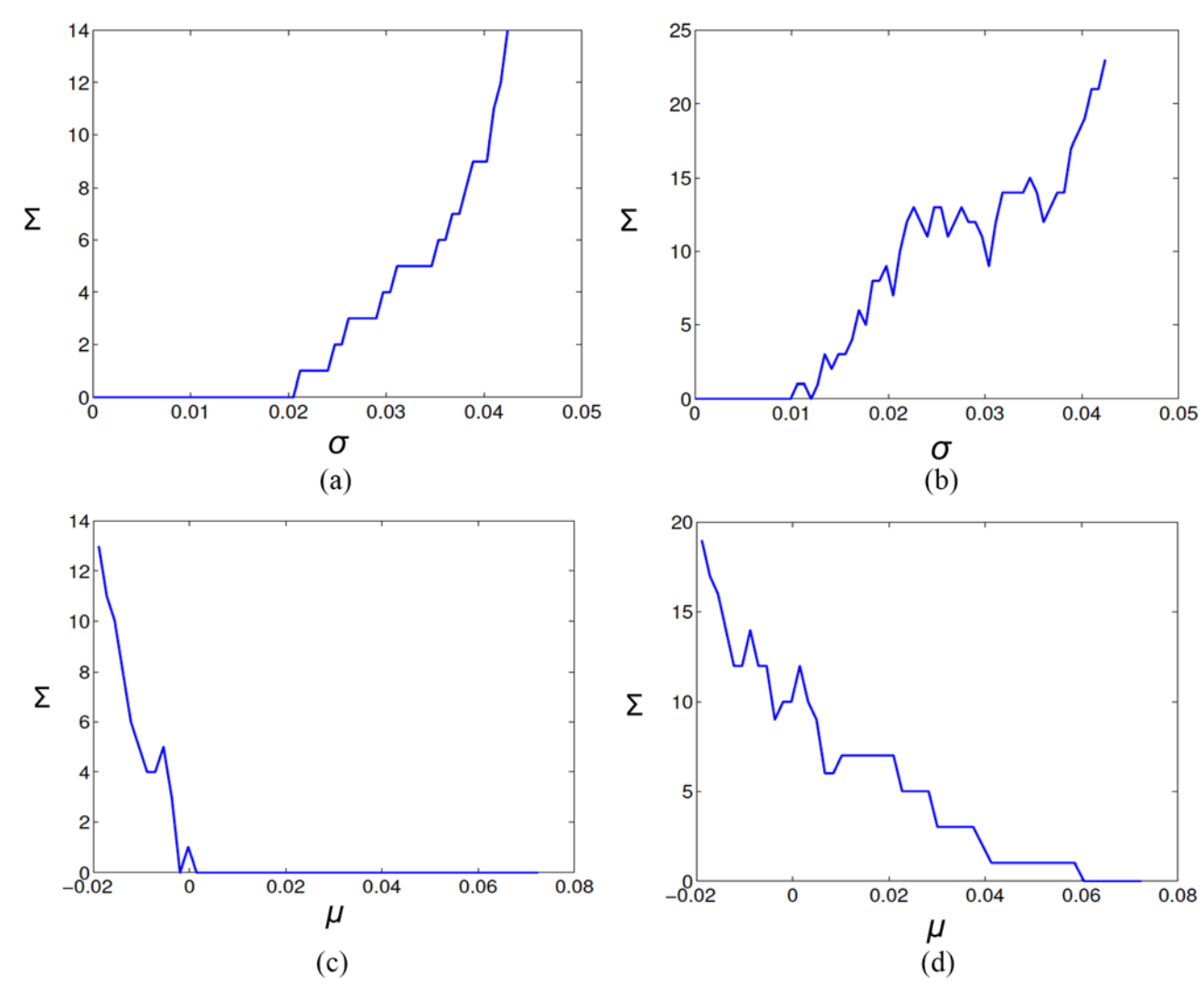

FIG. 12. Total number of outbreaks for a random sample path over the interval [0, 200]. Here $\Sigma$ represents the cardinality of the set $\{N = 0\}$. (a) $\mu = 0.0508$. (b) $\mu = 0.0015$. (c) $\sigma = 0.0099$. (d) $\sigma = 0.0219$.

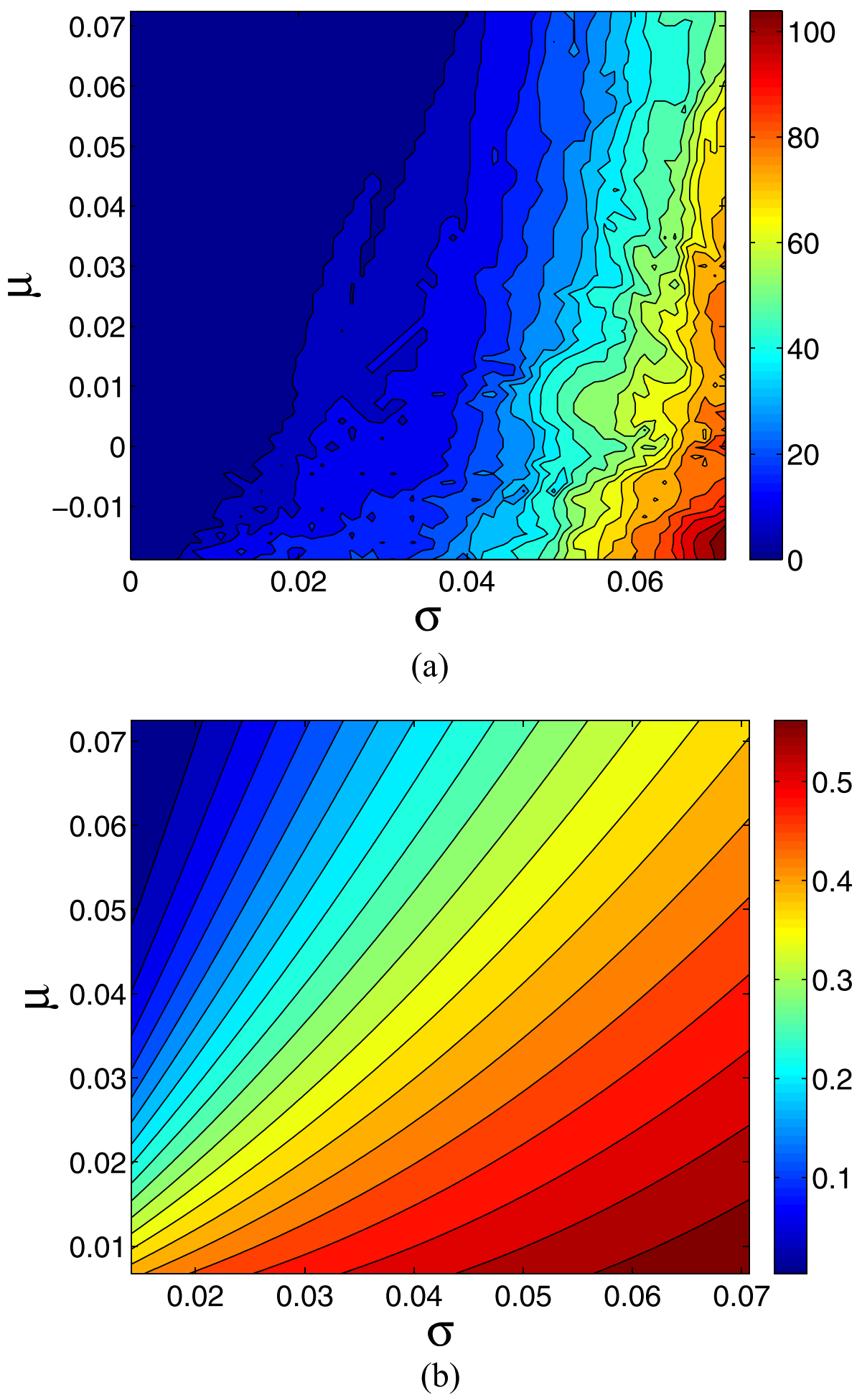

FIG. 13. (a) Total number of outbreaks for a random sample path over [0, 200]. (b) An approximate probability of repeated outbreaks. (a) Contour plot of $\Sigma$. (b) Contour plot of $\Phi(-\kappa)$.

Larger values of $N$ indicate longer periods of small amplitude fluctuations in the densities of the species, while very small values of $N$ represent a higher probability of recurring outbreaks. The asymptotic geometric distribution of $N$ gets more prominent when the distance to the Hopf bifurcation, $\mu$, is very small, approximately of the order $\varepsilon$, or the noise intensity is moderately large. However, when $\varepsilon < \mu < \sqrt{\varepsilon}$ as in Figure 5(a), where $\varepsilon$ is chosen to be 0.05, the distribution of $N$ seems to qualitatively resemble the distribution of return time of larch budmoths in the subalpine forest as studied in Ref. 12. With a proper choice of the other ecological parameters, it may be possible to replicate population dynamics of forest pests or of small mammals very closely. We plan to do a statistical analysis to estimate the parameters in our future work.

In the regime of weak noise, we note that $\mathbb{E}(N)$ is large, meaning that the sample paths, on an average, make large numbers of SAOs before exhibiting an LAO. In this case, the population outbreaks (if any) are separated by longer time periods. This situation is favorable in an ecosystem as the long-lived SAOs persist, thereby letting population fluctuations stay well under control. For larger noise amplitude, $\mathbb{E}(N)$ is low which suggests possibilities of frequent outbreaks. To control outbreaks, monitoring environmental fluctuations can be helpful. Environmental monitoring programs can be used to assess the quality of an environment and record the changing environmental parameters. These data can be then used to correlate with the fluctuations in the birth rate or the death rate of the species.

Our analysis can be used to find an estimate of the probability of two consecutive outbreaks (see Figure 13(b)). Negative values of $\kappa$ indicate that this probability is at least half, a situation which must be avoided for ecological reasons. Keeping the timescale ratio $\varepsilon$ and the noise amplitude $\sigma$ fixed, while treating $\mu$ as the control parameter, we observe that recurrences of outbreaks decrease as $\mu$ increases as shown in Figure 12. The rate of decrease depends on the noise amplitude. On the other hand, keeping $\varepsilon$ and $\mu$ fixed, while varying $\sigma$, we observe that the number of spikes increases with $\sigma$ (see Figure 12). A contour plot of recurrences of outbreaks as a function of $\mu$ and $\sigma$ is presented in Figure 13(a). When $\mu < 0$ but of order $\varepsilon$, the deterministic model exhibits a canard explosion.[28] Depending on the magnitude of $\mu$, in this regime in the presence of white noise, we expect to observe a wide variety of behaviors such as noise induced MMOs, bursting relaxation oscillations, coherence resonance, etc. Similar dynamics were observed in the stochastic FitzHugh–Nagumo model in Ref. 24 with a constant diffusion matrix. We are currently investigating this regime.

Understanding the mechanisms behind high amplitude population fluctuations and taking measures to control epidemic outbreaks have always formed an important field of study in population biology. Empirical studies have provided strong evidence that large scale fluctuations in population densities of small mammals and insects are found in environments with strong seasonal or between-year variations. The effect of seasonal cycles in a stochastic environment remains a subject for future investigation. It is interesting to note that the simple two species stochastic model considered in this paper can reproduce the population fluctuations of small mammals or insects as observed in Nature. The theoretical analysis presented in this work may provide useful information in identifying ecologically sensitive parameter regimes and taking measures to prevent outbreaks accordingly.

## ACKNOWLEDGMENTS

The author would like to thank the anonymous reviewers and Dr. S. Chakraborty Thakur from the University of California at San Diego for their suggestions that helped in improving the manuscript. The author would also like to thank the Faculty Scholarship Support Program at Georgia College for supporting this research.

## APPENDIX A: ASYMPTOTIC STABILITY OF THE FIXED POINT P OF SYSTEM (5)

From the definition of $\alpha^\star$ and under condition (12), it follows that

$$0 < 1-\beta-2\alpha^\star < \varepsilon(1-d). \tag{A1}$$

Hence by (A1), we obtain that

$$\begin{aligned}\det(J)|_{\alpha^\star} &= \frac{\alpha^\star(1-\alpha^\star)}{\varepsilon}\left[\frac{\varepsilon\beta(1-\alpha^\star)-\alpha^\star(1-\beta-2\alpha^\star)^2}{\varepsilon(1-\alpha^\star)(\beta+\alpha^\star)^2}\right]\\ &> \frac{\alpha^\star(1-\alpha^\star)}{\varepsilon}\left[\frac{\beta(1-\alpha^\star)-\varepsilon\alpha^\star(1-d)^2}{(1-\alpha^\star)(\beta+\alpha^\star)^2}\right]\\ &> \frac{\alpha^\star(1-\alpha^\star)}{\varepsilon}\left[\frac{\beta(1-2\alpha^\star)}{(1-\alpha^\star)(\beta+\alpha^\star)^2}\right]\\ &> \frac{\alpha^\star\beta^2}{\varepsilon(\beta+\alpha^\star)^2} > 0.\end{aligned}$$

Thus, the fixed point $P$ undergoes a Hopf bifurcation when $\alpha=\alpha^\star$. For $h>h^\star$, we will show that $P$ is asymptotically stable provided that (13) holds. Note that (13) implies (12). For $\alpha>\alpha^\star$, we have from (A1) that

$$\begin{aligned}&2\alpha\alpha^\star+2(\alpha+\alpha^\star)\beta-\beta+\varepsilon\beta+\beta^2,\\ &>2\alpha^{\star 2}+4\alpha^\star\beta-\beta+\varepsilon\beta+\beta^2\\ &>2\alpha^{\star 2}+\beta(1-\varepsilon-\beta)+2\varepsilon\beta d>0,\end{aligned}$$

for sufficiently small $\varepsilon$. Hence it follows from (8) that $Tr(J)<0$ for all such $\alpha$.

We note that the $y$-nullcline $hy=\frac{x}{\beta+x}-d$ is a decreasing function of $h$; hence, for all $h$ satisfying $h^\star<h<\tilde{h}$, we have that $\alpha^\star<\alpha<(1-\beta)/2$ and hence it follows from (9) that

$$\det(J) > \frac{\alpha(1-\alpha)}{\varepsilon}\left[\frac{4\beta}{(1+\beta)^2}-\tilde{h}\varepsilon(1-d)\right]>0,$$

provided $\varepsilon$ is sufficiently small. On the other hand, when $\alpha>(1-\beta)/2$, then it clearly follows from (9) that $\det(J)>0$. Thus $P$ is asymptotically stable for $h>h^\star$ as long as (13) holds.

## APPENDIX B: NORMAL FORM OF SYSTEM (5)

Here we derive the normal form of system (5) by carrying out the following transformations:

(1) We first rescale time as $dt\mapsto(\beta+x)dt$ to obtain an orbitally equivalent system

$$\begin{cases}\varepsilon\dot{x}=x\big((1-x)(\beta+x)-y\big)\\ \dot{y}=y\big(x-(d+hy)(\beta+x)\big).\end{cases} \tag{A2}$$

(2) We next perform an affine transformation $x=u+\alpha^\star$, $y=v+(1-\alpha^\star)(\beta+\alpha^\star)$ that translates the bifurcation point to the origin to obtain

$$\begin{cases}\varepsilon\dot{u}=-\alpha^\star v+\varepsilon(\alpha^\star-d(\beta+\alpha^\star))u+(1-\beta-3\alpha^\star)u^2\\ \quad -uv-u^3\\ \dot{v}=(1-\alpha^\star)(\beta+\alpha^\star)\delta\\ \quad +(1-\alpha^\star)(\beta+\alpha^\star)\big(1-d-h(1-\alpha^\star)(\beta+\alpha^\star)\big)u\\ \quad +\Big(\alpha^\star-d(\beta+\alpha^\star)-2h(1-\alpha^\star)(\beta+\alpha^\star)^2\Big)v\\ \quad +\big(1-d-2h(1-\alpha^\star)(\beta+\alpha^\star)\big)uv-h(\beta+\alpha^\star)v^2\\ \quad -huv^2,\end{cases}$$

where $\delta$ is defined by Equation (16). Using the relation in (6), (16) can be rewritten as

$$\begin{aligned}\delta=-(\beta+\alpha^\star)(\alpha-\alpha^\star)\Bigg[&\frac{\beta}{(\beta+\alpha)(\beta+\alpha^\star)}\\ &-h\big(1-\beta-(\alpha+\alpha^\star)\big)\Bigg],\end{aligned} \tag{A3}$$

which is of the order $(\alpha-\alpha^\star)$ near the bifurcation point, and thus is related to the distance from the Hopf point. In fact, it can be shown that $\delta$ is negative when $h>h^\star$, and positive when $h<h^\star$, provided that $\varepsilon>0$ is small and that (13) holds.

(3) We next scale space and time by $u=\sqrt{\varepsilon}p$, $v=-\varepsilon q$ and $t=\sqrt{\varepsilon}t'$ to obtain

$$\begin{cases}\dot{p}=\alpha^\star q+(1-\beta-3\alpha^\star)p^2+\sqrt{\varepsilon}[(\alpha^\star-d(\beta+\alpha^\star))p\\ \quad +pq-p^3]\\ \dot{q}=-(1-\alpha^\star)(\beta+\alpha^\star)\dfrac{\delta}{\sqrt{\varepsilon}}-\gamma p-(h(1-\alpha^\star)(\beta+\alpha^\star)^2\\ \quad -\delta)\sqrt{\varepsilon}q+\big(1-d-2h(1-\alpha^\star)(\beta+\alpha^\star)\big)\varepsilon pq\\ \quad +h(\beta+\alpha^\star)\varepsilon^{3/2}q^2+h\varepsilon^2pq^2,\end{cases}$$

where $\dot{}=d/dt'$ and $\gamma$ is defined by (17).

(4) We then perform a nonlinear transformation $\alpha^\star q=\xi-(1-\beta-3\alpha^\star)p^2+\frac{\alpha^\star\gamma}{2(1-\beta-3\alpha^\star)}$, which has the effect of straightening out the $\dot{p}=0$ nullcline, resulting into

$$\begin{cases}\dot{p}=\xi+\dfrac{\alpha^\star\gamma}{2(1-\beta-3\alpha^\star)}+\sqrt{\varepsilon}\Bigg[\Bigg(\alpha^\star-d(\beta+\alpha^\star)\\ \quad +\dfrac{\gamma}{2(1-\beta-3\alpha^\star)}+\dfrac{\xi}{\alpha^\star}\Bigg)p-\dfrac{(1-\beta-2\alpha^\star)}{\alpha^\star}p^3\Bigg]\\ \dot{\xi}=-\alpha^\star(1-\alpha^\star)(\beta+\alpha^\star)\dfrac{\delta}{\sqrt{\varepsilon}}+2(1-\beta-3\alpha^\star)p\xi\\ \quad +\sqrt{\varepsilon}\Bigg[\Big(\delta-h(1-\alpha^\star)(\beta+\alpha^\star)^2\Big)\Bigg(\xi+\dfrac{\alpha^\star\gamma}{2(1-\beta-3\alpha^\star)}\\ \quad -(1-\beta-3\alpha^\star)p^2\Bigg)+2(1-\beta-3\alpha^\star)\Bigg(-d(\beta+\alpha^\star)\\ \quad +\alpha^\star+\dfrac{\gamma}{2(1-\beta-3\alpha^\star)}+\dfrac{\xi}{\alpha^\star}\Bigg)p^2\\ \quad -\dfrac{2}{\alpha^\star}(1-\beta-2\alpha^\star)(1-\beta-3\alpha^\star)p^4\Bigg]\\ \quad +\varepsilon\Bigg[\big(1-d-2h(1-\alpha^\star)(\beta+\alpha^\star)\big)p\Bigg(\xi+\dfrac{\alpha^\star\gamma}{2(1-\beta-3\alpha^\star)}\\ \quad -(1-\beta-3\alpha^\star)p^2\Bigg)\Bigg]+\varepsilon^{3/2}h(\beta+\alpha^\star)\Bigg(\xi+\dfrac{\alpha^\star\gamma}{2(1-\beta-3\alpha^\star)}\\ \quad -(1-\beta-3\alpha^\star)p^2\Bigg)^2+\varepsilon^2\dfrac{h}{\alpha^\star}p\Bigg(\xi+\dfrac{\alpha^\star\gamma}{2(1-\beta-3\alpha^\star)}\\ \quad -(1-\beta-3\alpha^\star)p^2\Bigg)^2.\end{cases}$$

(5) Finally, we rescale $p$ and $q$ by defining new variables $l$ and $z$ as $l = \frac{1-\beta-3\alpha^\star}{\alpha^\star\gamma}p$, $z = -\frac{1-\beta-3\alpha^\star}{\alpha^\star\gamma}\xi$, which yields system (14). The variables $A_i$, $i = 0,\ldots 4$ are defined by

$$A_0(l,z) = \left(\alpha^\star - d(\beta+\alpha^\star) + \frac{\gamma}{2(1-\beta-3\alpha^\star)} - \frac{\gamma z}{(1-\beta-3\alpha^\star)}\right)l - \frac{\alpha^{\star 2}\gamma^2(1-\beta-2\alpha^\star)}{(1-\beta-3\alpha^\star)^2}l^3, \tag{A4}$$

$$A_1(l,z) = \left(\delta - h(1-\alpha^\star)(\beta+\alpha^\star)^2\right)\left(z - \frac{1}{2} + \alpha^\star\gamma l^2\right) - 2\alpha^\star\gamma l^2\left(\alpha^\star - d(\beta+\alpha^\star) + \frac{\gamma}{2(1-\beta-3\alpha^\star)} - \frac{\gamma z}{1-\beta-3\alpha^\star}\right) + \frac{2\alpha^{\star 2}\gamma^3(1-\beta-2\alpha^\star)l^4}{(1-\beta-3\alpha^\star)^2}, \tag{A5}$$

$$A_2(l,z) = \frac{\alpha^\star\gamma l}{1-\beta-3\alpha^\star}\left(1 - d - 2h(1-\alpha^\star)(\beta+\alpha^\star)\right)\left(z - \frac{1}{2} + \alpha^\star\gamma l^2\right), \tag{A6}$$

$$A_3(l,z) = \frac{-h(\beta+\alpha^\star)\gamma}{1-\beta-3\alpha^\star}\left(z - \frac{1}{2} + \alpha^\star\gamma l^2\right)^2, \tag{A7}$$

$$A_4(l,z) = \frac{-h\alpha^\star\gamma^2 l}{(1-\beta-3\alpha^\star)^2}\left(z - \frac{1}{2} + \alpha^\star\gamma l^2\right)^2. \tag{A8}$$

## APPENIDX C: NORMAL FORM OF SYSTEM (4)

**Proposition IX.1.** *In the new variables $(L, Z)$ and on the new time scale, system (4) takes the form*

$$\begin{cases} dL = \left[\frac{1}{2} - Z + B_0(L,Z)\sqrt{\varepsilon}\right]dt' + \hat{\sigma_1}C(L)G_1(L)dB_1(t') \\ dZ = \left[\hat{\mu} + 2\alpha^\star\gamma LZ + B_1(L,Z)\sqrt{\varepsilon} + B_2(L,Z)\varepsilon + B_3(L,Z)\varepsilon^{3/2} + B_4(L,Z)\varepsilon^2\right]dt' \\ \quad + C(L)\left[\hat{\sigma_1}G_2(L)dB_1(t') + \hat{\sigma_2}G_3(L,Z)dB_2(t')\right], \end{cases} \tag{A9}$$

*where*

$$\hat{\sigma_1} = \frac{1-\beta-3\alpha^\star}{\gamma\varepsilon^{3/4}}\sigma_1,\ \hat{\sigma_2} = \frac{1-\beta-3\alpha^\star}{\gamma\varepsilon^{3/4}}\sigma_2, \tag{A10}$$

$$\hat{\mu} = \mu - \alpha^\star\gamma(\beta+\alpha^\star)\hat{\sigma_1}^2, \tag{A11}$$

$$C(L) = \sqrt{\beta + \alpha^\star + \frac{\sqrt{\varepsilon}\alpha^\star\gamma L}{1-\beta-3\alpha^\star}}, \tag{A12}$$

$$B_0(L,Z) = A_0(L,Z), \tag{A13}$$

$$B_1(L,Z) = A_1(L,Z) - \frac{\alpha^\star\gamma^2\hat{\sigma_1}^2(\beta+\alpha^\star+\alpha^\star L)}{1-\beta-3\alpha^\star}, \tag{A14}$$

$$B_2(L,Z) = A_2(L,Z) - \frac{\alpha^\star\gamma^3\hat{\sigma_1}^2(\beta+3\alpha^\star)L^2}{(1-\beta-3\alpha^\star)^2}, \tag{A15}$$

$$B_3(L,Z) = A_3(L,Z) - \frac{\alpha^{\star 2}\gamma^4\hat{\sigma_1}^2L^3}{(1-\beta-3\alpha^\star)^3}, \tag{A16}$$

$$B_4(L,Z) = A_4(L,Z), \tag{A17}$$

$$G_1(L) = \left(1 + \frac{\sqrt{\varepsilon}\gamma L}{1-\beta-3\alpha^\star}\right), \tag{A18}$$

$$G_2(L) = -2\alpha^\star\gamma\hat{\sigma_1}\left(1 + \frac{\sqrt{\varepsilon}\gamma L}{1-\beta-3\alpha^\star}\right), \tag{A19}$$

$$G_3(L,Z) = \hat{\sigma_2}\left((1-\alpha^\star)(\beta+\alpha^\star) + \frac{\varepsilon\gamma\left(Z - \frac{1}{2} + \alpha^\star\gamma L^2\right)}{1-\beta-3\alpha^\star}\right), \tag{A20}$$

*and $A_i$, $i = 0, 1, \ldots, 4$ are defined by (A4)–(A8).*

## APPENDIX D: PROOF OF PROPOSITION VI.1

Solving Equation (20) by variation of constants yields

$$Z_T = Z_0 + e^{\frac{\gamma\alpha^\star T^2}{2}}\left(\hat{\mu}\int_{t_0}^{T} e^{-\frac{\gamma\alpha^\star s^2}{2}}ds - \gamma\alpha^\star\hat{\sigma_1}\sqrt{\beta+\alpha^\star}\int_{t_0}^{T} se^{-\frac{\gamma\alpha^\star s^2}{2}}dB_1(s) + (1-\alpha^\star)(\beta+\alpha^\star)\hat{\sigma_2}\sqrt{\beta+\alpha^\star}\int_{t_0}^{T} se^{-\frac{\gamma\alpha^\star s^2}{2}}dB_2(s)\right). \tag{A21}$$

From the definition of $P$, it follows that $e^{2\gamma\alpha^\star P^2} = (c_0\hat{\mu})^{-a}$. Let $T = 4P$ so that when $t' = T$, $L_T \approx P$. We are interested in the location of $Z_T$. If $Z_T < b < 0$ for some large $|b|$, then the system will perform a large excursion which will correspond to a spike. On the other hand, if $Z_T > 0$, then the system will exhibit an SAO because the term $\hat{\mu}t'$ pushes the sample paths upwards causing them to revolve around $P$. We note from (A21) that the random variable $Z_T$ is Gaussian with expectation and variance given by

$$\mathbb{E}\{Z_T\} = Z_0 + \hat{\mu}e^{2\gamma\alpha^\star P^2}\int_{-2P}^{2P} e^{-\frac{\gamma\alpha^\star s^2}{2}}ds, \tag{A22}$$

$$Var(Z_T) = (\beta+\alpha^\star)e^{\gamma\alpha^\star P^2}\left(\int_{-2P}^{2P} s^2e^{-\gamma\alpha^\star s^2}ds\right)\left(\gamma^2\alpha^{\star 2}\hat{\sigma_1}^2 + (1-\alpha^\star)^2(\beta+\alpha^\star)^2\hat{\sigma_2}^2\right), \tag{A23}$$

respectively. Thus for any real $H$

$$\mathbb{P}\{Z_T \le -H\} = \int_{-\infty}^{H}\frac{e^{-\frac{(Z-\mathbb{E}\{Z_T\})^2}{2Var(Z_T)}}}{\sqrt{2\pi Var(Z_T)}} = \Phi\left(-\frac{H + \mathbb{E}\{Z_T\}}{\sqrt{Var(Z_T)}}\right). \tag{A24}$$

Note from (A22) and (A23) and the definition of $P$ that

$$\mathbb{E}\{Z_T\} \approx Z_0 + \hat{\mu} e^{2\gamma\alpha^\star P^2}\sqrt{\frac{2\pi}{\gamma\alpha^\star}}$$
$$= Z_0 + \hat{\mu}(c_0\hat{\mu})^{-a}\sqrt{\frac{2\pi}{\gamma\alpha^\star}}, \quad \text{(A25)}$$

$$Var(Z_T) \le (\beta+\alpha^\star)e^{\gamma\alpha^\star P^2}\sqrt{\frac{\pi}{\gamma\alpha^\star}} \times \left(\frac{\gamma\alpha^\star\hat{\sigma_1}^2}{2} + (1-\alpha^\star)^2(\beta+\alpha^\star)^2\hat{\sigma_2}^2\right)$$
$$= (\beta+\alpha^\star)(c_0\hat{\mu})^{-2a}\sqrt{\frac{\pi}{\gamma\alpha^\star}} \times \left(\frac{\gamma\alpha^\star\hat{\sigma_1}^2}{2} + (1-\alpha^\star)^2(\beta+\alpha^\star)^2\hat{\sigma_2}^2\right). \quad \text{(A26)}$$

Applying the bounds given by (A25) and (A26) on (A24), the result follows.